\documentclass[12pt]{amsart}

\usepackage{epsfig,color}

\usepackage{amssymb}

\newtheorem{thm}{Theorem}[section]
\newtheorem{prop}[thm]{Proposition}
\newtheorem{defin}[thm]{Definition}
\newtheorem{corol}[thm]{Corollary}

\newtheorem{lem}[thm]{Lemma}

\newtheorem{rem}[thm]{Remark}

\newtheorem{exa}[thm]{Example}

\newtheorem*{theorem1}{Theorem 1}
\newtheorem*{theorem2}{Theorem 2}

\newtheorem*{proposition}{Proposition 1}
\newtheorem*{conjecture1}{Conjecture 1}
\newtheorem*{conjecture2}{Conjecture 2}
\newtheorem*{corollary}{Corollary 1}

\newcommand{\FF}{{\mathcal F}}

\newcommand{\G}{{\mathcal G}}

\newcommand{\mm}{{\mathcal M}}

\newcommand{\N}{{\mathbb N}}

\newcommand{\R}{{\mathbb R}}
\newcommand{\RR}{{\mathbb R}^2}

\newcommand{\Z}{{\mathbb Z}}
\newcommand{\ZZ}{\mathbb Z^2}

\newcommand{\bo}{\partial} %%%%%%%%%% BOUNDARY

\newcommand{\card}{\mbox{card}} %%%%%%%%%% CARDINALITY
\newcommand{\id}{\text{Id}}   %%%%%%%%%%%%%%%% IDENTITY
\newcommand{\inn}{\text{in}}        %%%%%%%%%%%%%%%%% INNER

\newcommand{\len}{\text{length}}     %%%%%%%%%%%%%%%%% LENGTH
\newcommand{\out}{\text{out}}   %%%%%%%%%%%%%%%%% OUTER

\newcommand{\per}{\text{per}}  %%%%%%%%%%%%% PERIODIC

\newcommand{\sgn}{\text{sgn}}
\newcommand{\stab}{\text{Stab}}

\newcommand{\tot}{\text{tot}}    %%%%%%%%%%%%%%%%  TOTAL

\newcommand{\ep}{\varepsilon}

\newcommand{\tal}{\tilde{\al}}

\newcommand{\tc}{\tilde{c}}
\newcommand{\td}{\tilde{d}}

\newcommand{\tg}{\tilde{g}}\newcommand{\tga}{\tilde{\ga}}

\newcommand{\al}{\alpha}
\newcommand{\be}{\beta}
\newcommand{\ga}{\gamma}\newcommand{\Ga}{\Gamma}

\newcommand{\ve}{\varepsilon}
\newcommand{\ip}{\varphi}
\newcommand{\la}{\lambda}
\newcommand{\si}{\sigma}

\def\qed{\hfill\vrule height 5pt width 5 pt depth 0pt}

\date{\today}

\begin{document}

\title[Blocking and Flatness]
{Insecurity for compact surfaces of positive genus}

\author{Victor Bangert and Eugene Gutkin}

\address{Mathematisches Institut\\
Albert-Ludwigs-Universit\"at\\
Eckerstrasse 1\\
79104 Freiburg im Breisgau\\
Germany}

\email{bangert@email.mathematik.uni-freiburg.de}

\address{Nicolaus Copernicus University (UMK)\\ 
and Mathematics Institute of the
Polish Academy of Sciences (IMPAN)\\
Chopina 12/18\\87-100 Torun\\
Poland}

\email{gutkin@mat.umk.pl,\ gutkin@impan.pl}

\keywords{compact riemannian surface, minimal geodesic, closed geodesic, 
free homotopy class, blocking set, insecurity, flat torus, surface of genus greater than one}\

\subjclass{53C22,57M10,37E40}

%\thanks{K.B.~is partially supported by N.S.F.~grant DMS-0408704.}

\maketitle

\tableofcontents

%\date{March 26, 2005}
%\date{\today}

\begin{center}
\begin{minipage}{120mm}

\vskip 0.1in
\baselineskip 0.2in

{{\bf Abstract.} A pair of points in a riemannian manifold $M$ is
secure if the geodesics between the points can be blocked by a
finite number of point obstacles; otherwise the pair of points is insecure. 
A manifold is secure  if all pairs of points in $M$
are secure. A manifold is insecure if there exists an insecure point pair,  
and totally insecure if all point pairs are insecure.

Compact, flat manifolds are secure. A standing conjecture says
that these are the only secure, compact riemannian manifolds.  We
prove this for surfaces of genus greater than zero. We also prove
that a closed surface of genus greater than one with any riemannian
metric and a closed surface of genus one with generic metric are
totally insecure.

}
%\medskip\noindent

\end{minipage}
\end{center}

\section{Introduction}  \label{intro}
We begin by describing our setting and establishing the
terminology. By a riemannian manifold $(M,g)$ we will mean a
complete, connected, infinitely differentiable riemannian
manifold. We will view geodesics in $(M,g)$ as curves, $c:I\to M$,
parameterized by arclength, where $I\subset\R$ is any nontrivial
interval. The set $c(I)\subset M$ is the trace of $c$. If
$I=[a,b]$, we will also say that $c(I)$ is a geodesic segment; the
points $x=c(a),y=c(b)$ are the endpoints of the geodesic. If $z\in
M$ is an interior point of $c(I)$, we say that $c$ {\em passes
through $z$}.

For any pair $x,y\in M$ (including $y=x$) let $G(x,y)$ be the set
of geodesic segments in $(M,g)$ with endpoints $x,y$. We say that
these geodesics {\em join} $x$ with $y$. A finite set $B\subset
M\setminus\{x,y\}$ is a {\em blocking set} for the pair $x,y$ if
every geodesic in $G(x,y)$ passes through a point $b\in B$. We
will also say that $B$ {\em blocks $x$ and $y$ away from each
other}.

A pair $x,y\in M$ is {\em insecure} if these points cannot be
blocked away from each other. Otherwise, the pair $(x,y)$ is
secure. A riemannian manifold is insecure if it contains an
insecure pair of points. Thus, $(M,g)$ is {\em secure}  if any
point in it can be blocked away from any point, including itself.
Finally, $(M,g)$ is {\em totally insecure} if every pair $x,y\in
M$  is insecure. See \cite{Gut05} for a motivation of this
terminology.

Which compact riemannian manifolds are secure? The only examples so far
are the flat manifolds \cite{GS06}. Researchers in the subject
believe in the following statement \cite{BG08,LS07}.
\begin{conjecture1}   \label{main_conj}
A compact riemannian manifold is secure if and only if it is flat.
\end{conjecture1}
Restricting the setting, we state a counterpart of
Conjecture~\ref{main_conj} for tori.
\begin{conjecture2}   \label{tori_conj}
A riemannian torus is secure if and only if it is flat.
\end{conjecture2}
In this note we obtain several results concerning security of
closed riemannian surfaces. In particular, we establish
Conjecture~\ref{tori_conj} for two-dimensional tori. See
Theorem~~\ref{torus_thm}.

\begin{theorem1}   \label{tori_thm}
A two-dimensional riemannian torus is secure
if and only if it is flat.
\end{theorem1}

\vspace{3mm}

Security and insecurity are preserved under finite coverings
\cite{GS06}. This observation allows us to restrict the discussion
to orientable manifolds. Our next result concerns the insecurity
of closed surfaces of genus greater than one. See
Theorem~~\ref{high_genus_thm}.

\begin{theorem2}   \label{total_thm}
A closed riemannian surface of  genus greater than one is totally
insecure.
\end{theorem2}

Since surfaces of genus  greater than one do not admit flat
riemannian metrics, this Theorem provides evidence for
Conjecture~\ref{main_conj}.

\medskip

An insecure manifold, in particular, a nonflat two-torus, may have
secure pairs of points. In section~~\ref{revolut} we analyze the
security of point pairs on two-dimensional tori of revolution. This
class of riemannian tori contains the round euclidean tori of revolution
$T^2\subset\R^3$. Let $E_{\inn}\subset T^2$ be the
inner equator. Our Corollary~~\ref{round_tor_cor} says the following.

\begin{proposition}   \label{tor_rev_prop}
Let $T^2\subset\R^3$ be a round euclidean torus of revolution. Let $x,y\in T^2$ be any pair of
points. Then it is secure if and only if $x,y\in E_{\inn}$.
\end{proposition}
\medskip

Tori of revolution are special. By our
Corollary~~\ref{open_dense_cor}, a generic riemannian two-torus
is totally insecure. Here is the precise statement.
\begin{corollary}    \label{generic_cor}
The set of riemannian two-tori contains a $C^2$-open and
$C^{\infty}$-dense subset of totally insecure tori.
\end{corollary}
\vspace{3mm}

%We will now briefly discuss the literature.
Conjecture~\ref{main_conj} holds for locally symmetric spaces
\cite{GS06}. Moreover, compact locally symmetric spaces of
noncompact type are totally insecure \cite{GS06}. This is true, in
particular, for compact surfaces of constant negative curvature. A
geometric argument showing this is sketched in
\cite{GS06}. Let $M$ be a surface as above, let $x,y\in M$ be a
pair of points, and let $C\subset M\setminus\{x,y\}$ be a periodic
geodesic. The argument in question constructs an infinite sequence
of geodesics $\ga_n\in G(x,y)$. As $n\to\infty$, the geodesics
$\ga_n$, their lengths going to infinity, spend almost all their
time in ever smaller vicinity of $C$, never intersecting it. Thus,
through any point $z\in M\setminus\{x,y\}$ passes at most a finite
number of geodesics $\ga_n$.

\medskip

Our proofs of Theorem~~\ref{tori_thm} and Theorem~~\ref{general}
use a similar idea, although the two approaches differ
considerably in detail. Since we are working with arbitrary
riemannian metrics, we cannot use the hyperbolicity of the
geodesic flow, which was crucial for the analysis in \cite{GS06}. 
Instead, we use the classical
results of Morse \cite{Mo} and Hedlund \cite{He}, as well as more
recent results \cite{Ba1}, \cite{Ba2}, \cite{In}, on minimal
geodesics in surfaces to construct infinite sequences of joining
geodesics with a similar behavior. See section~\ref{constr}
and section~\ref{key}. Our
Proposition~~\ref{suffice_prop} says that once we have an infinite
sequence of geodesics joining points $x,y$ and satisfying certain
conditions, the pair $(x,y)$ is insecure.  The strategy of proving
Theorem~~\ref{tori_thm} and Theorem~~\ref{general} is to construct
in each case sequences of geodesics satisfying the requirements of
Proposition~~\ref{suffice_prop}.

\medskip

Another approach to insecurity for riemannian manifolds is based
on connections between the insecurity and the growth rate of
the number of joining geodesics. See \cite{BG08}, \cite{LS07}, and \cite{Gut09}.
This approach works well under additional assumptions, e. g., that
the manifold $(M,g)$ has no conjugate points. Then we can use the
well known relationships between the growth rate of the number
$n_T(x,y)$ of joining geodesics having length $\le T$ and the
growth of $\pi_1(M)$, as well as those between $n_T(x,y)$ and the
topological  entropy of $(M,g)$.

%\vspace{6mm}

\medskip

\noindent{\em Acknowledgements.} The work of V. Bangert was partially supported by SFB/Transregio 71 
``Geometrische Partielle Differentialgleichungen''. E. Gutkin did some of the work 
while visiting FIM at ETH, Zurich, the Department of Mathematics at UCLA, and
the Mathematisches Institut at Albert-Ludwigs-Universit\"at in Freiburg im Breisgau.
It is a pleasure to thank these institutions for hospitality and financial support.

\medskip

%The exposition  is organized as follows. {\bf Later}
%\noindent{\bf Acknowledgements.} A substantial portion of the present work took place during
%the  second author's visit to the Institute for Mathematical Research at ETH in Zurich.
%It is a pleasure to thank the Institute and its director, Marc Burger, for hospitality and
%financial support.

%\newpage

\section{Rays, corays, and Busemann functions} \label{general}
In this section we review basic material on rays, corays and
Busemann functions in riemannian manifolds. We will use the
notation $(M^n,g)$ for riemannian manifolds, suppressing $g$
and/or $n$ whenever this causes no confusion. We denote by
$d(\cdot,\cdot)$ the riemannian distance on $M$. We will view
geodesics as parameterized curves $c(t),\,t\in I$, where
$I\subset\R$ is a nontrivial, possibly infinite interval, and $t$
is an arclength parameter. We will call the set $c(I)\subset M$
the {\em trace of $c$}.
%(See section~\ref{key}, and especially the paragraph
%preceding the proof of Lemma~\ref{key_lem} for precise notational
%conventions.)

\medskip
\begin{defin}   \label{ray_def}
\vspace{1mm}
Let $I\subset\R$ be an interval and let $c:I\rightarrow M$ be a
geodesic.
\vspace{1mm}
\begin{itemize}
\item[(a)]
The geodesic $c$ is  {\em minimal} if $d(c(t), c(s))=|t-s|$ for
all $s,t\in I$.
\item[(b)]
A {\em ray} is a minimal geodesic $c:\R_{\ge 0}\rightarrow M$.
\item[(c)]
Let $c:\R_{\ge 0}\rightarrow M$ be a ray, and let $C\subset M$ be its trace.

\noindent A ray $\tilde{c}$ is a {\em coray} to $c$ if there exists a
sequence of minimal geodesics $c_n:[0, L_n]\rightarrow M$ with
$\lim_{n\rightarrow\infty} L_n=\infty$, such that
$\lim_{n\rightarrow\infty} \dot{c}_n(0)=\dot{\tilde{c}}(0)$ and
$c_n(L_n)\in C$ for all $n\in\N$.
\end{itemize}
\end{defin}

\medskip

Throughout the paper we will denote by $L(\ga)$ the length of a
curve.
\begin{defin}   \label{homot_min_def}
A geodesic $c:I\rightarrow M$ is {\em homotopically minimal} if
its lifts to the universal riemannian covering of $M$ are minimal.
\end{defin}

Equivalently, a geodesic $c:I\rightarrow M$ is homotopically
minimal iff the following holds. Let $s,t\in I$ and let $s<t$. Let
$\ga:[s,t]\to M$ be a curve satisfying $\ga(s)=c(s),\ga(t)=c(t)$
and such that the curves $\ga$ and $c|_{[s,t]}$ are homotopic with
fixed endpoints. Then
$$
L(c([s,t]))\le L(\ga).
$$
\medskip

Throughout the paper we will use the following basic result. It is
essentially due to M. Morse. See \cite{Mo}, Theorems~~9 and~~10.

\begin{thm}              \label{minim_length_thm}
Let $M$ be an orientable surface, and let $c:\R\to M$ be a
periodic geodesic with period $L>0$. If the closed curve
$c|_{[0,L]}$ has minimal length among all closed curves freely
homotopic to $c|_{[0,L]}$, then $c$ is homotopically minimal.
\end{thm}

\medskip

Taking limits of minimal geodesics of finite length, we obtain the
following basic fact.

\medskip\noindent
\begin{prop}  \label{ray_lem}
A complete riemannian manifold $(M,g)$ carries a ray if and only
if it is not compact. If $c$ is a ray in $(M,g)$ and $p\in M$,
then there exists a coray $\tilde{c}$ to $c$ with
$\tilde{c}(0)=p$.
\end{prop}
\medskip
%\noindent
%We will now recall the notion of {\bf Busemann functions.}
%
\begin{defin}   \label{busem_def}
Let $c:\R_+\to M$ be a ray. Its {\em Busemann function},
$B_c:M\to\R$, is defined by
\begin{equation}   \label{busem_eq}
B_c(p) = \lim_{t\rightarrow\infty}\left[d(p, c(t))-t\right].
\end{equation}
\end{defin}

By the triangle inequality, the function $t\rightarrow d(p,
c(t))-t$ is monotonically decreasing.\footnote{In general, not
strictly.} Also by the triangle inequality, it satisfies
$-d(c(0),p)\le d(p, c(t))-t.$ Thus, the limit in
equation~\eqref{busem_eq} exists.

\medskip\noindent
\begin{lem}  \label{busem_lip_lem}
Let $c:\R_+\rightarrow M$ be a ray, and let $p, q\in M$ be
arbitrary points. Then
$$
|B_c(p)-B_c(q)|\le d(p, q).
$$
\begin{proof}
Apply the triangle inequality to the triangle with corners $p,
q,c(t)$, and take the limit $t\to\infty$.
\end{proof}
\end{lem}

By Lemma~\ref{busem_lip_lem}, any Busemann function is lipschitz,
with the lipschitz constant $1$.

\medskip\noindent
\begin{prop}   \label{coray_thm}
Let $c:\R_+\to M$ be a ray. A geodesic $\tc:\R_+\to M$ is a coray
to $c$ if and only if for all $s,t\in\R_+$ the equation
$$
B_c(\tc(t))-B_c(\tc(s))=s-t
$$
holds.
\begin{proof}
This follows from equations (22.16) and (22.20) in \cite{Bu}.
\end{proof}
\end{prop}

\medskip
%\noindent
We use Proposition~\ref{coray_thm} to relax the requirements in
Definition~\ref{ray_def}.

\medskip
%\noindent
%
\begin{lem}  \label{coray_lem}
Let $c:\R_+\to M$ be a ray, and let $C\subset M$ be its trace. Let
$L_n\to\infty$ be a positive sequence. Let $c_n:[0,L_n]\to M$ be
minimal geodesics such that $\lim_{n\to\infty} d(c_n(L_n),C)=0$.

If $\tc:\R_+\to M$ is a geodesic such that $\dot{\tc}(0)$ is a
point of accumulation of the sequence $\dot{c}_n(0)$, then $\tc$
is a coray to $c$.
\begin{proof}
We assume without loss of generality that $\lim_{n\to\infty}
\dot{c}_n(0)=\dot{\tilde{c}}(0)$. For every $n\in\N$ there is a
number $t_n\in[0,\infty)$ such that the sequence
$\ve_n=d(c_n(L_n),c(t_n))$ converges to zero. Since the geodesics
$c_n$ are minimal, the condition $\lim_{n\to\infty}L_n=\infty$
implies that $\lim_{n\to\infty}t_n=\infty$. By the triangle
inequality, for all $n$ and any $t\in[0,L_n]$ we have
\begin{equation}                  \label{coray_eq}
L_n - t - \ve_n\le d(c_n(t), c(t_n))\le L_n - t + \ve_n.
\end{equation}
Let $s, t>0$ be arbitrary. Using that $\lim\dot{c}_n(0) =
\dot{\tilde{c}}(0)$ and equation~\eqref{coray_eq}, we have
\begin{eqnarray*}        %\label{coray_EQ}
B_c(\tilde{c}(t)) - B_c(\tilde{c}(s))   =
\lim_{n\to\infty} [d(\tilde{c}(t), c(t_n)) - d(\tilde{c}(s), c(t_n))]\\
=\lim_{n\to\infty} [d(c_n(t), c(t_n)) - d(c_n(s),
c(t_n))]\ge\lim_{n\to\infty} (s-t-2\ve_n)=s-t.
\end{eqnarray*}
Combining this inequality with Lemma~\ref{busem_lip_lem}, we
obtain
\begin{equation}        \label{coray_EQ}
B_c(\tilde{c}(t)) - B_c(\tilde{c}(s))=s-t.
\end{equation}
The claim now follows from Proposition~\ref{coray_thm}.
\end{proof}
\end{lem}
\section{Outline of the proof that nonflat two-tori are insecure}   \label{idea}
For the benefit of the reader, we will outline the main ideas in
the proof of Theorem~\ref{torus_thm} which is Theorem~~\ref{intro}
in the Introduction. The methods used in this
proof are also typical for several other statements in this paper.
Let $(T^2,\overline{g})$ be a non-flat two-dimensional torus; our
goal is to find a pair of points in $(T^2,\overline{g})$ that
cannot be blocked away from each other by a finite blocking set.

%Let $\overline{g}$ be a
%riemannian metric on the $2$-torus $T^2=\RR/\ZZ$.

\medskip

By a classical theorem of E.~Hopf \cite{Ho}, a riemannian two-torus is flat if and only if
it has no conjugate points. Thus, the
torus $(T^2, \overline{g})$ has conjugate points. Then, by a
theorem of N.~Innami, there exists
a nontrivial free homotopy class $\alpha$ of closed curves such that
$(T^2, \overline{g})$ cannot be foliated by geodesics in  $\alpha$.
See \cite{In}, Corollary 3.2; see also the proof of Theorem 6.1 in \cite{Ba2}.

\medskip

Let $\mm_{\al}^{\per}$ be the set of periodic geodesics of minimal length in the
class $\alpha$. By results that go back to M.~Morse \cite{Mo} and
G.~Hedlund \cite{He}, these geodesics do
not self-intersect and are pairwise disjoint. Generically,
$\mm_{\al}^{\per}$ consists of a single geodesic.

\medskip

%{\bf Note: The bad notation for sets in the torus and the covering plane
%originates from the text below. Here is what I propose. Not to use ``bar'' for
%objects in the torus. For instance, use $(T^2, g)$ instead of $(T^2, \overline{g})$.
%Use ``tilde'' for the liftings to the plane. For instance, use
%$(\RR, \tg)$ instead of $(\RR, g)$.
%Instead of $Y$ for the cylinder in $(T^2, g)$ use $A$. Then set
%$\tA=S$ for the lift to $\RR$. Thus, it will be $\mm_{\al}$ in the torus, and
%$\tilde{\mm}_{\al}$ in the plane. Same proposal applies to the geodesics in $S$
%discussed below in the paragraph preceding Figure 1.}

The geodesics in $\mm_{\al}^{\per}$ foliate a compact, proper subset,
$N\subset T^2$.  Let $Z\subset T^2$ be a connected component of
$T^2\setminus N$; let $p,q\in Z$ be any pair of points. We will
show that the pair $p,q$ is insecure, i. e., that we cannot block
$p$ away from $q$ by a finite blocking set.

\medskip

We denote by $(\RR,g)$ the riemannian universal covering; let
$\pi:(\RR,g)\to(T^2,\overline{g})$ be the projection.
%We will use
%the representation $T^2=\RR/\ZZ$.
Let $S\subset\RR$ be a connected component of $\pi^{-1}(Z)$. Then
$S$ is an open strip. The boundary $\bo S$ is a disjoint union of
traces of two geodesics, $c_0:\R\to(\RR,g)$ and
$c_1:\R\to(\RR,g)$. By Theorem~~\ref{minim_length_thm}, the
geodesics $c_0,c_1$ are minimal. Let $C_0,C_1$ be the respective
traces; then $\bo S=C_0\cup C_1$.

\medskip

Let $P,Q_0\in S$ be arbitrary points such that $\pi(P)=p,\pi(Q_0)=q$.
Using the action of the stabilizer of $S$ in $\pi_1(T^2)=\ZZ$, we produce an infinite sequence of points
$Q_1,\dots,Q_n,\ldots\in S$ such that  $\pi(Q_n)=q$ and the sequence of distances $L_n=d(P,Q_n)$
goes to infinity.
Let now $\tc_n:[0,L_n]\to S$ be a sequence of minimal geodesics such that
$\tc_n(0)=P$ and $\tc_n(L_n)=Q_n$.

\begin{figure}[htbp]
\begin{center}
\input{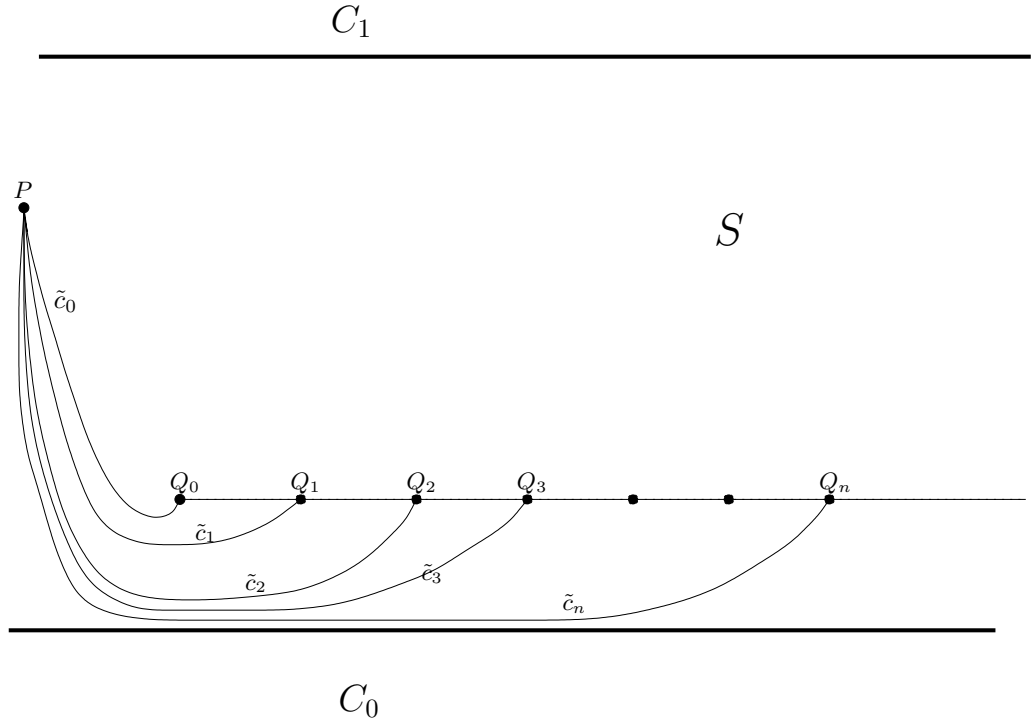}
\caption{A sequence of minimal geodesics in the universal covering
hose projections to the torus cannot be blocked by a finite point
set. }
\label{fig1}
\end{center}
\end{figure}

Lemma~\ref{key_lem} in section~\ref{key}
implies that most of the
time the geodesics $\tc_n$ are close to $\bo S$. More precisely,
for any $\ve>0$ there exists $T=T(\ve)>0$ such that
for all $t\in[T,L_n-T]$ the points $\tc_n(t)$ are $\ve$-close to $\bo S$.
Figure~\ref{fig1} illustrates the behavior of this sequence of geodesics.

\medskip

Set $c_n=\pi\circ\tc_n$. Then the geodesics $c_n:[0,L_n]\to Z$
join the points $p,q$. Let $z\in Z\setminus\{p,q\}$ be an arbitrary point.
The preceding discussion implies that at most a finite
number of the geodesics $c_n$ passes through $z$. On the other hand, if $z\in T^2\setminus Z$,
no geodesic in our sequence passes through it.
Thus, any point $z\in T^2$
can block at most a finite number of joining geodesics
in the infinite sequence $c_n$.
Hence, we cannot block the points $p,q$ away from each other by
a finite blocking set.

We will now illustrate the preceding discussion with the example
of round euclidean tori of revolution.
\begin{exa}   \label{tor_rev_exa}
{\em
Let $0<r<R$, and set $C=C(r,R)=\{(x,0,z):(x-R)^2+z^2=r^2\}$.
This is a circle of radius $r$ in the $xz$-plane. The round euclidean torus of revolution, $T^2=T^2(r,R)\subset\R^3$,
is obtained by revolving $C$ about the $z$-axis. Points in $C$ yield the 
{\em circles of latitude} in $T^2$. Let $\al$ be their free homotopy class.

Let the inner equator $E_{\inn}\subset T^2$ (resp. the outer equator $E_{\out}\subset T^2$)
be the  circle of latitude generated by the point
$(R-r,0,0)\in C$ (resp. $(R+r,0,0)\in C$). Their lengths are $2\pi(R-r)$ and $2\pi(R+r)$ respectively.
These are the only circles of latitude that are geodesics in the round euclidean torus.

Specializing the preceding discussion to the torus of revolution, we have
$N=E_{\inn}$ and $Z=T^2\setminus E_{\inn}$. The set $\mm_{\al}^{\per}$ consists of a single geodesic; 
although the tori of revolution are very special, this is the generic situation
for two-dimensional riemannian tori. See section~~\ref{generic}. By the preceding argument,
no points $p,q\in T^2\setminus E_{\inn}$ can be blocked
away from each other by a finite blocking set.

In section~~\ref{revolut} we will completely analyze the blocking for
a class of two-dimensional tori containing the above tori of revolution.

}
\end{exa}

\vspace{4mm}

\section{Minimal geodesics in an admissible strip}   \label{constr}
We will use the notation of section~\ref{idea}; in particular, we
use the identification $(T^2, \overline{g})=(\RR,g)/\ZZ$. If
$S\subset\RR$, we denote by
$$\stab(S) = \{j\in\ZZ:\,S+j=S\}
$$
the stabilizer of $S$.
%More precisely,
%$\stab(S)=\{j\in\ZZ:\,S+j=S\}$.
Recall that a nonzero vector $j\in\ZZ$ is {\em prime} if there do
not exist $n\in\N$, $n\ge 2$, and $k\in\ZZ$ such that $j=nk$.

\medskip

We will use the notation $\bar{X}$ for the closure of a set.

\medskip

%\medskip\noindent
%
\begin{prop}     \label{strip_prop}
Let $(T^2, \overline{g})$ be a nonflat riemannian torus. Let
$(\RR, g)$ be its universal covering, and let $\pi:(\RR,
g)\to(T^2, \overline{g})$ be the projection.

Then there exists a connected open set $S\subset\RR$ with
totally geodesic boundary, such that the following statements
hold.
\begin{itemize}
\item[(a)]
The group $\stab(S)$ is generated by a prime vector.
%$k\in\ZZ\setminus\{0\}$, such that $S+k=S$; equivalently, $S$ is
%invariant under translation by $k$.
\item[(b)]
If $j\in\ZZ\setminus\stab(S)$, then
$(S+j)\cap\overline{S}=\emptyset$.
\item[(c)]
The boundary of $S$ has two connected components, say $C_0$ and
$C_1$. There are minimal geodesics $c_0, c_1:\R\to(\RR,g)$ whose
traces are $C_0$ and $C_1$, respectively.
\item[(d)]
Let $c:\R\to S$ be a geodesic such that $\pi\circ c:\R\to T^2$ is
periodic. Then $c$ is not minimal.
\end{itemize}
\begin{proof}
By Corollary 3.2 in \cite{In}, there exists a nontrivial free
homotopy class, say $\alpha$, of closed curves in $T^2$ having the
following property: There does not exist a family of closed
geodesics in the class $\alpha$ whose traces foliate $T^2$. We can
assume that $\alpha$ is prime.

%We identify $H_1(T^2)$ with $\ZZ$; let $k\in\ZZ$ be the vector
%corresponding to $\al$. By Theorem 6.6 in \cite{Ba1}, we can
%choose $\al$ so that $k$ is a prime vector.

%\medskip

Let $L$ be the minimal length of a curve in $\al$; we denote by
$\mm_{\al}^{\per}$ the set of closed geodesics in the class $\al$ having
length $L$. Clearly, $\mm_{\al}^{\per}\ne\emptyset$. By Theorem 6.5 and
Theorem 6.6 in \cite{Ba1}, the trace of every $c\in\mm_{\al}^{\per}$ is
an embedded curve in $T^2$. Moreover, if $c,\tc\in\mm_{\al}^{\per}$, then
either their traces are disjoint or $c$ and $\tc$ coincide up to a
translation of the parameter.

\medskip

Let $N$ be the union of the traces of geodesics in $\mm_{\al}^{\per}$. By
our choice of $\al$, the set $N\subset T^2$ is a proper, nonempty, closed
subset. Let $Z$ be a connected component of $T^2\setminus N$. Let
$\bo Z=\overline{Z}\setminus Z$ be its boundary. Then either $\bo
Z$ is the trace of a geodesic in $\mm_{\al}^{\per}$ or $\bo Z$ is the
union of traces of two geodesics in  $\mm_{\al}^{\per}$.\footnote{We
point out that $Z$ is homeomorphic to the cylinder
$S^1\times(0,1)$.}

Let $S$ be a connected component of $\pi^{-1}(Z)\subset\RR$. Then
the boundary of $S$ is the union of the traces of two geodesics
$c_0,c_1:\R\to\RR$ such that $\pi\circ c_0$ and $\pi\circ c_1$
belong to $\mm_{\al}^{\per}$. Let $k\in\ZZ$ correspond to $\al$. Then for
all $t\in\R$ we have
\begin{equation}   \label{period_eq}
c_0(t+L)=c_0(t)+k,\,c_1(t+L)=c_1(t)+k.
\end{equation}
%
%\medskip

Theorem~~\ref{minim_length_thm} implies that $c_0$ and $c_1$ are
homotopically minimal geodesics. The remaining statements in (a), (b), and (c)
now follow by elementary topological arguments; claim (d) follows
from Theorem 6.7 in \cite{Ba1}.
\end{proof}
\end{prop}

\medskip

In what follows an open set $S\subseteq\RR$ satisfying the
conditions stated in Proposition~~\ref{strip_prop} will be called
an {\em admissible strip}; the projection $Z=\pi(S)\subseteq T^2$
is an {\em admissible cylinder}.

\section{The key lemma}   \label{key}
We will use the setting and the notation of section~\ref{constr}.
In particular, $S$ will denote an admissible strip. The following
statement is crucial in our proof of Theorem~\ref{torus_thm}. We
will refer to it as the {\em Key Lemma}.

\begin{lem}     \label{key_lem}
For any $\ve>0$ there exists $T=T(\ve)>0$ such that the following
holds. If $c:[0,L]\rightarrow S$ is a minimal geodesic and
$d(c(0),\partial S)\ge\ve$ then $d(c(t),\partial S)\le\ve$ for all
$t\in[T,L-T]$.
\end{lem}

%\noindent{\bf the way it is stated, $T=T(\eta,\delta)$, i. e., it
%does not depend on a particular $c:[0,L]\to S$. Is this right?}
%\medskip

The proof of Lemma~\ref{key_lem} is based on the results of
M.~Morse \cite{Mo} about minimal geodesics in $S$ and on a result
from \cite{Ba2} concerning the rays in $S$.
%We will deduce
%Lemma~\ref{key_lem} from these results by compactness arguments.
%Our proofs are indirect: We will assume that claims fail, and
%obtain contradictions.
We need a few technical lemmas.

\begin{lem}     \label{dist_lem}
Let  $c:[0,\infty)\rightarrow S$ be a ray. Then $\lim_{t\to\infty}
d(c(t),\bo S)=0$.
\begin{proof}
The claim follows from Theorem 3.7 in \cite{Ba2}, interpreted as a
statement about minimal geodesics in $(\RR,g)$. See Example (1) on page 51
in \cite{Ba2} for details.
\end{proof}
\end{lem}
\medskip

Throughout this section we will use the following notational conventions.
With any geodesic $c:\R\to\overline{S}$
we will associate two geodesics $c_{\pm}:\R_+\to\overline{S}$ as follows.
The geodesic $c_+$ is the restriction of $c$ to the positive half-line.
We define the geodesic $c_-$ by $c_-(t)=c(-t)$.

%If $a:I\to\overline{S}$ is any geodesic,
%we denote by $A\subset\overline{S}$ the corresponding geometric curve. Formally,
%$A=\{a(t):\,t\in I\}$. Applying this convention to $c:\R\to\overline{S}$

We will denote by $C,C_+,C_-\subset\overline{S}$ the respective traces of
$c,c_+,c_-$.

\vspace{6mm}

\begin{defin}     \label{coherent_def}
Let $c_0,c_1:\R\to\overline{S}$ be two geodesics such that their
traces $C_0,C_1$ are the two components of $\bo S$. We say that
the geodesics $c_0,c_1$ are {\em coherently oriented} if for any
time sequence $t_n\to\infty$ the two point sequences
$c_0(t_n),c_1(t_n)\in\bo S$ converge to the same end of
$\overline{S}$.
\end{defin}

\begin{figure}[htbp]
\begin{center}
\input{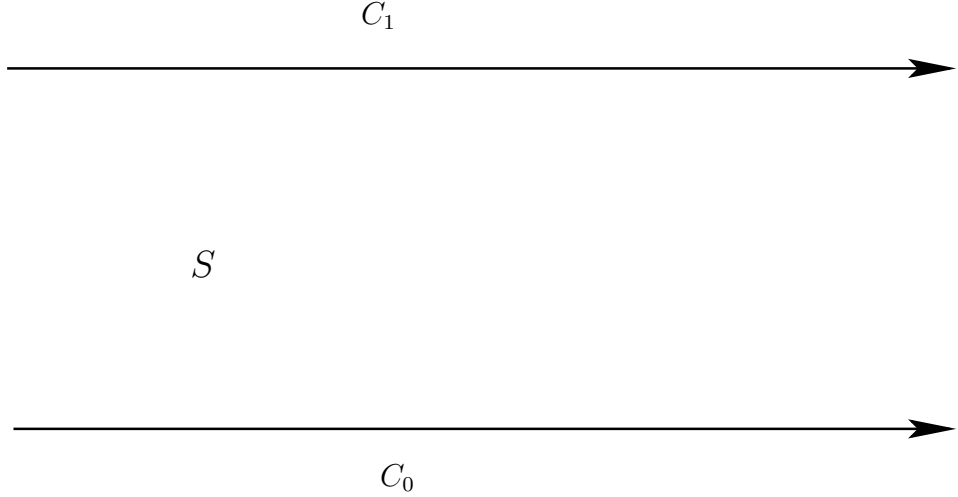}
\caption{A strip with coherently oriented boundary components.}
\label{fig2}
\end{center}
\end{figure}

Figure~\ref{fig2} illustrates Definition~\ref{coherent_def}.
We will also say that $c_0,c_1:\R\to\overline{S}$ are {\em coherent parameterizations}
of $\bo S$.

\medskip

\begin{lem}     \label{component_lem}
Let $\bo S =C_0\cup C_1$, where $c_0,c_1:\R\to\overline{S}$ are
coherently oriented. Let the geodesics $c_{0,\pm}:\R_+\to\bo S$
and $c_{1,\pm}:\R_+\to\bo S$ be as above; let $C_{0,\pm}\subset\bo
S$ and $C_{1,\pm}\subset\bo S$ be the respective traces.

Let now $c:\R\to S$ be a minimal geodesic.
Then, switching $c_0$ with $c_1$ and reversing the orientation of $c$, if need be, we have
\begin{equation}    \label{ends_eq}
\lim_{t\to-\infty}d(c(t),C_{0,-})=0,\ \lim_{t\to+\infty}d(c(t),C_{1,+})=0.
\end{equation}
\begin{proof}
By Lemma~\ref{dist_lem}, $c(t)$ converges to $\bo S$ as $|t|\to\infty$.
By Theorem 15 in \cite{Mo} or by Theorem 6.7 in \cite{Ba1}, the equation
$\lim_{t\to-\infty}d(c(t),C_{0,-})=0$ implies $\lim_{t\to+\infty}d(c(t),C_{1,+})=0$.
\end{proof}
\end{lem}
%
%
%{\bf I propose to make the following remark}
%
\begin{rem}    \label{connect_rem}
{\em Let $M$ be a riemannian manifold. If $c:I\to M$ is a geodesic, its {\em inverse} is the
geodesic $c^{-1}:-I\to M$ defined by $c^{-1}(t)=c(-t)$. Lemma~\ref{component_lem} is
equivalent to the following geometric fact.

\noindent Let $Z\subset T^2$ be an admissible cylinder. Assume,
for simplicity of exposition, that the closure of $Z$ is a proper
subset of $T^2$. Let
$\overline{c}_0,\overline{c}_1:\R\to(T^2,\overline{g})$ be the
periodic geodesics in the homotopy class $\al$ whose respective
traces are the two components of the boundary $\bo Z$.

Let now $c:\R\to Z$ be a geodesic whose lift $\tc:\R\to S$ is minimal.
Then  $c$ is  a heteroclinic connection either between
$\overline{c}_0$ and $\overline{c}_1$ or
between $\overline{c}_0^{-1}$ and $\overline{c}_1^{-1}$.

}
\end{rem}

%
%\begin{rem}     \label{dir_rem}
%{\em Since $\lim_{t\to\infty} d(c(-t),c(t))=\infty$, the rays
%$c_+$ and $c_-$ are asymptotic to $C_0$ and $C_1$ respectively
%``in opposite directions''. }
%\end{rem}
%
\medskip

%From now on we assume that the boundary components of $S$ are coherently oriented,
%and the conclusions of Lemma~\ref{component_lem} are satisfied.
Our next lemma says that if a ray in $S$ is a coray  to a ray
in a boundary component of $S$, then it is asymptotic to this component.

\begin{figure}[htbp]
\begin{center}
\input{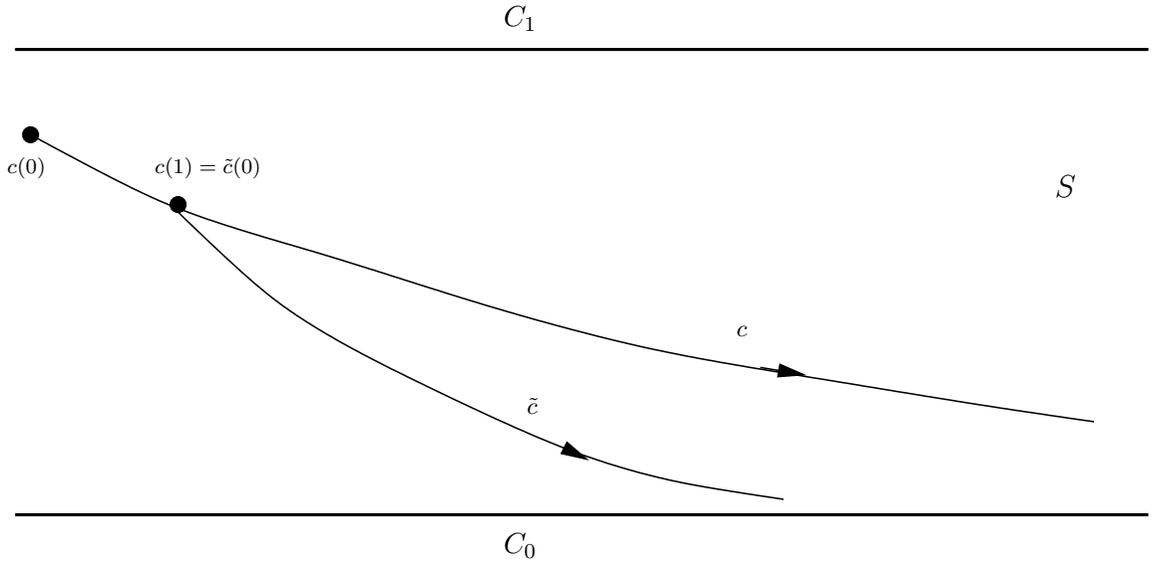}
\caption{Illustration to the proof of Lemma~\ref{asympt_lem}: A
configuration that does not exist.} \label{fig3}
\end{center}
\end{figure}

\begin{lem}     \label{asympt_lem}
Let $c_0:\R\to\overline{S}$ be a geodesic whose trace is one of the components of  $\bo S$.

Let $c:\R_+\to S$ be a coray to $c_{0,+}$. Then $\lim_{t\to\infty}d(c(t),C_{0,+})=0$.
\begin{proof}
By Theorem 3.7 in \cite{Ba2}, for any $q\in S$ there exists a ray $\tc:\R_+\to S$ such that
$\tc(0)=q$ and $\lim_{t\to\infty}d(\tc(t),C_{0,+})=0$. Hence,
by Lemma~\ref{dist_lem} and Lemma~\ref{coray_lem},
$\tc$ is a coray to $c_{0,+}$.

Set $q=c(1)$, and let $\tc:\R_+\to S$ be as above. Thus, both $c$
and $\tc$ are corays to $c_{0,+}$; by construction, $\tc(0)=c(1)$.
The geodesic $t\mapsto c(1+t)$ is also a coray  to $c_{0,+}$ starting at
$q=c(1)=\tc(0)$. By Theorem 22.19 in \cite{Bu} or, by Corollary
3.8 in \cite{Ba2}, there is only one coray to $c_{0,+}$ starting
at $q$. Figure~\ref{fig3} shows a hypothetical configuration of
the rays $c$ and $\tc$ which cannot materialize.

Therefore, the ray $\tc$ satisfies
$\tc(t)=c(1+t)$. Since $\tc$ is asymptotic to $C_{0,+}$, the claim follows.
\end{proof}
\end{lem}
%

%\begin{center}
%Figure 3\\

%Illustration to Lemma~\ref{asympt_lem}: A configuration that cannot happen
%\end{center}

\medskip

We will now prove a preliminary variant of the Key Lemma.

\begin{figure}[htbp]
\begin{center}
\input{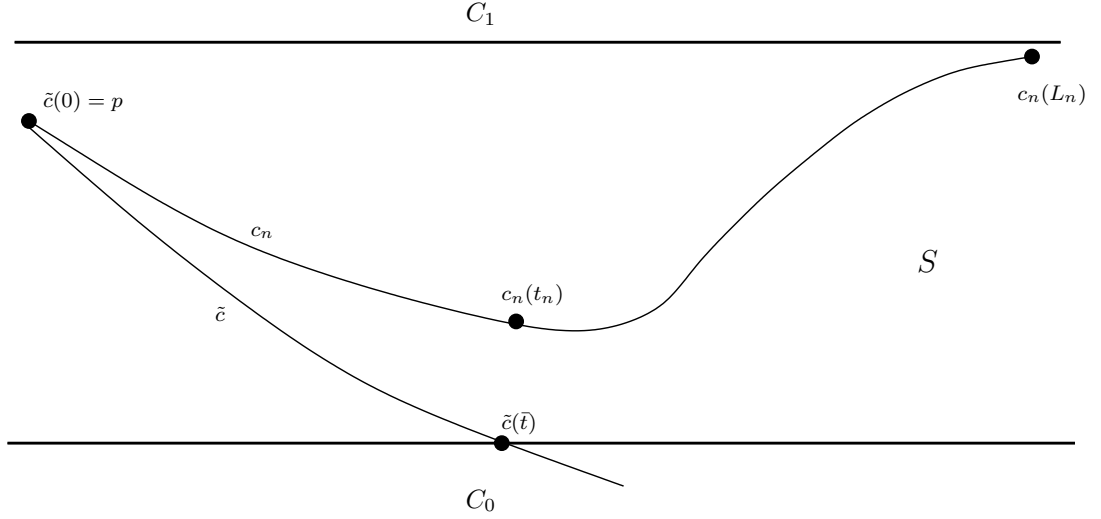}
\caption{Illustration to the proof of Lemma~\ref{prekey_lem}: The
behavior of the geodesic $\tc$ deduced from the assumption
$\lim_{n\to\infty}t_n<\infty$.} \label{fig4}
\end{center}
\end{figure}

\begin{lem}     \label{prekey_lem}
For any $\ve>0$ there exists $\eta=\eta(\ve)>0$ such that the
following holds. Let $0<L<\infty$, and let $c:[0,L]\to S$ be a
minimal geodesic such that $d(c(0),\partial S)\ge\ve$ and
$d(c(L),C_1)<\eta$. Then $d(c(t),C_0)\ge\eta$ for all $t\in[0,L]$.
\begin{proof}
Suppose that the claim fails. Then there exists $\delta>0$, a
sequence of minimal geodesics $c_n:[0,L_n]\to S$, and a sequence
$t_n\in[0,L_n]$ such that $d(c_n(0),\partial S)\ge\delta$,
$\lim_{n\to\infty}d(c_n(L_n),C_1)=0$, and
$\lim_{n\to\infty}d(c_n(t_n),C_0)=0$.

The closed strip $\overline{S}$ is invariant under the group $\stab(S)\simeq k\Z$
that acts on $\overline{S}$ by isometries. We have denoted this action by $z\mapsto z+rk$.
We will use the same notation for the corresponding action of  $\stab(S)$ on geodesics
in $S$. Then for any
integers $r_1,\dots,r_n,\ldots\in\Z$ the sequence of geodesics $\tc_n=c_n+r_nk$ satisfies
the above conditions. In view of this observation, and the compactness of the
quotient $\overline{S}/k\Z$, we assume without loss of generality
that the vectors $\dot{c}_n(0)$ converge to a limit vector, $v\in T^1(S,g)$;
let $p\in S$ be its footpoint.

We will now prove that $\lim_{n\to\infty}t_n=\infty$. If this
fails, then, by passing to a subsequence, if need be, we have
$\lim t_n=\overline{t}<\infty$. Let $\tc:\R\to\RR$ be the geodesic
with the initial vector $v$. Then $\tc(0)=p\in S$, and
$\tc(\overline{t})=q=\lim_{n\to\infty} c_n(t_n)\in C_0\subset\bo S$.
%We have $\tc(0)=p\in S$ and $\tc(\overline{t})=q\in C_0\subset\bo S$.
Since $\bo S$ is geodesic,
$\tc$ intersects it transversally at $q$. Thus, for $t>\overline{t}$ and
sufficiently close to $\overline{t}$, we have $\tc(t)\notin\overline{S}$.
Figure~\ref{fig4} illustrates the analysis.

On the other hand, $\liminf L_n\ge\overline{t}+d(C_0,C_1)$ implies that
$\tc(t)\in\overline{S}$ for all $t\in[0,\overline{t}+d(C_0,C_1)]$. In view of this
contradiction, $\lim t_n=\infty$.

\medskip

Since $\lim d(c_n(t_n),C_0)=0$ and $\lim t_n=\infty$, by
Lemma~\ref{coray_lem}, $\tc$ is a coray to $c_{0,+}$ or $c_{0,-}$.
Similarly,  $\lim d(c_n(L_n),C_1)=0$ and $\lim L_n=\infty$ imply,
by Lemma~\ref{coray_lem}, that $\tc$ is a coray to $c_{1,+}$ or
$c_{1,-}$. In view of Lemma~\ref{asympt_lem}, this is impossible.
\end{proof}
\end{lem}

\vspace{3mm}

We will now prove the Key Lemma. Recall that we view geodesics in
$\overline{S}$ as mappings $c:I\to\overline{S}$ of nontrivial
intervals $I\subset\R$. For $t\in I$ the velocity vectors
$\dot{c}(t)$ are unit tangent vectors in $(\RR,g)$. Thus,
$\len(c)=|I|$. If $0\in I$, we
will refer to $\dot{c}(0)$ as the {\em initial vector} of $c$.

\medskip
\begin{proof}({of Lemma~\ref{key_lem}})
Suppose the claim fails. Then for some $\ve>0$ there exists a
sequence of minimal geodesics $c_n:[0,L_n]\to S$ such that the
following conditions are satisfied:

\noindent i) For all $n\in\N$ we have $d(c_n(0),\bo S)\ge\ve$;

\noindent ii) For each $n$ there
is $t_n\in[n,L_n-n]$ so that $d(c_n(t_n),\bo S)>\ve$.

As in the proof of Lemma~\ref{prekey_lem}, we assume without loss
of generality that the velocity vectors $\dot{c}_n(t_n)$ converge
to a vector $v\in T^1(\overline{S},g)$. Let $\tc:\R\rightarrow S$
be the geodesic such that $v=\dot{\tc}(0)$. Since all of
$c_n:[0,L_n]\to S$ are minimal, and $t_n\in[n,L_n-n]$, we conclude
that $\tc:\R\rightarrow S$ is a minimal geodesic. By construction,
it satisfies $d(\tc(0),\bo S)\ge\ve$.

Let $\eta=\eta(\ve)>0$ be as in Lemma~\ref{prekey_lem}. By
Lemma~\ref{component_lem}, there are $s_0, s_1\in\R$ such that
$d(\tc(s_0),C_0)<\eta$ and $d(\tc(s_1),C_1)<\eta$.

Interchanging $C_0$ and $C_1$, if need be, we may assume that
$s_0<s_1$.
For any $t\in\R$ we have $\lim_{n\to\infty}c_n(t_n+t)=\tc(t)$.
In particular, $\tc(s_0)=\lim_{n\to\infty}c_n(t_n+s_0)$ and
$\tc(s_1)=\lim_{n\to\infty}c_n(t_n+s_1)$.
Therefore, for sufficiently large $n$  the inequalities $d(c_n(t_n+s_0), C_0)<\eta$
and $d(c_n(t_n+s_1), C_1)<\eta$ hold. Besides,
for sufficiently large $n$  we have $0<t_n+s_0<t_n+s_1<L_n$.

Let $n\in\N$ be any index such that the above conditions hold. Set
$L=t_n+s_1$, and let $c:[0,L]\to S$ be the restriction of $c_n$ to
$[0,t_n+s_1]$. Then $d(c(0),\bo S)\ge\ve$ and $d(c(L),C_1)<\eta$.
But we also have $d(c(t_n+s_0),C_0)<\eta$, and $t_n+s_0\in(0,L)$.
By Lemma~\ref{prekey_lem}, this is impossible.
\end{proof}
\medskip

\begin{figure}[htbp]
\begin{center}
\input{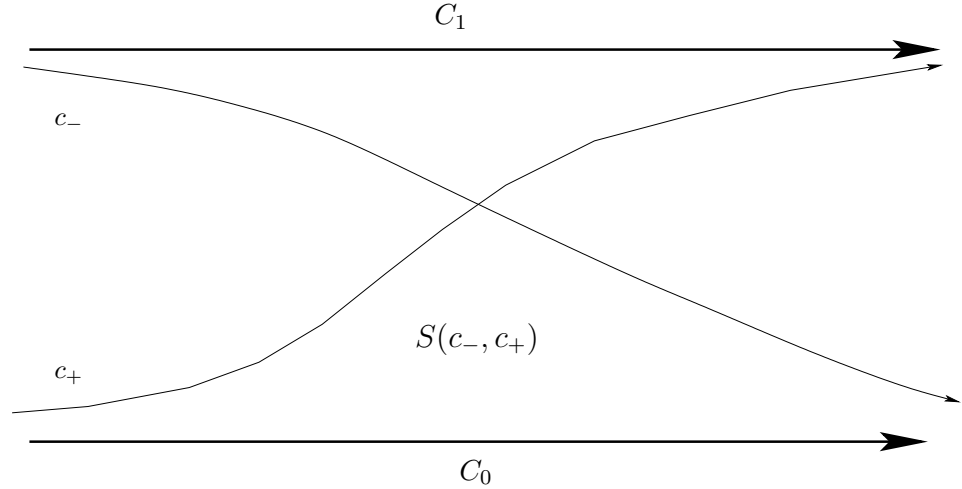}
\caption{Illustrates the proof of Lemma~~\ref{compl_lem}.}
\label{new1}
\end{center}
\end{figure}

We will also need the following complement to
Lemma~~\ref{key_lem}.
\begin{lem}     \label{compl_lem}
For any $\ve>0$ there exists $T=T(\ve)$ such that the following
holds. If $c:[0,L]\to\overline{S}$ is a minimal geodesic such that
$d(c(0),C_1)\ge\ve$ and $c(L)\in C_0$ then $d(c(t),C_0)\le\ve$ for
all $t\in[T,L]$.
\begin{proof}
By Theorem 14 in \cite{Mo} (see also Theorem 6.8 in \cite{Ba1}),
there exist minimal geodesics $c_{\pm}:\R\to S$ such that
$\lim_{t\to-\infty}d(c_+(t),C_{0,-})=0$ and
$\lim_{t\to\infty}d(c_+(t),C_{1,+})=0$, while
$\lim_{t\to-\infty}d(c_-(t),C_{1,-})=0$ and
$\lim_{t\to\infty}d(c_-(t),C_{0,+})=0$. Hence the geodesics $c_-$
and $c_+$ intersect; being minimal, they intersect exactly once.
Let $t_-,t_+\in\R$ be determined by $c_-(t_-)=c_+(t_+)$. We denote
by $S(c_-,c_+)\subset S$ the connected component of
$S\setminus(c_-(\R)\cup c_+(\R))$ whose boundary is $\bo
S(c_-,c_+)=C_0\cup c_+((-\infty,t_+])\cup c_-([t_-,\infty))$.
See Figure~\ref{new1}.

\medskip

If $k\in\stab(S)$ then $c_+-k$ and $c_-+k$ have the same
properties as $c_+$ and $c_-$. If $k$ is not the identity map,
then
$$
\cup_{r\in\Z}S(c_- + rk,c_+ - rk)\ =\ S.
$$
A straightforward compactness argument shows that for every
$\ve>0$ there exists $r\in\Z$ such that the following holds: If
$p\in S$ and $d(p,C_1)\ge\ve$ then there exists $j\in\stab(S)$
such that
$$
p+j\in S(c_- + rk,c_+ - rk).
$$

We set $c_-=c_-+rk,\,c_+=c_+-rk$ and use the preceding statement
for $p=c(0)$. Thus, translating $c$ by some $j\in\stab(S)$, if
necessary, we may assume that $c(0)\in S(c_-,c_+)$. By the
minimality of geodesics $c,c_-,c_+$ and the inclusions $c(0)\in
S(c_-,c_+),\,c(L)\in C_0$, we conclude that $c(t)\in S(c_-,c_+)$
for all $t\in [0,L]$. Now the conditions
$\lim_{t\to-\infty}d(c_+(t),C_{0,-})=0$ and
$\lim_{t\to\infty}d(c_-(t),C_{0,+})=0$ imply the claim. 
\end{proof}
\end{lem}

\medskip

\medskip

\begin{rem}       \label{compl_rem}
{\em The claim of Lemma~~\ref{compl_lem} is geometrically obvious, 
although a formal proof is somewhat
cumbersome. For the benefit of the reader, we indicate how to
formalize the last step in the preceding argument. Using the normal
exponential map of $C_0$, we represent the ends
$c_+((-\infty,-T_0])$ and $c_-([T_0,\infty))$ of $c_+$ and $c_-$
as graphs over $C_0$ for sufficiently large $T_0$. Then, for
sufficiently large $t\in[0,L]$ the point $c(t)$ is located between
$C_0$ and one of these graphs. 
%See Figure~\ref{new2}.

}
\end{rem}
%

%
%\begin{figure}[htbp]
%\begin{center}
%\input{newfig2.pstex_t}
%\caption{Another illustration for the proof of
%Lemma~~\ref{compl_lem}.} \label{new2}
%\end{center}
%\end{figure}
%

\medskip

\section{Nonflat two-dimensional tori are insecure}  \label{proof}
In this section we prove that every nonflat two-torus has an open
set of insecure point pairs. We begin with an auxiliary
proposition; it provides a general sufficient condition for a pair
of points in a riemannian manifold to be insecure.

\begin{prop}     \label{suffice_prop}
Let $p,q$ be points in a complete riemannian manifold $M$. Let
$A\subset M$ be a nonempty closed set. Suppose that there exists a
sequence of geodesics $c_n:[0,L_n]\to M$ from $p=c_n(0)$ to
$q=c_n(L_n)$ with the following properties:
\begin{itemize}
\item [a)]
We have $\lim_{n\to\infty}L_n=\infty$;
\medskip
\item [b)]
If $n\in\N$ and $t\in(0,L_n)$ then $c_n(t)\notin A$;
\medskip
\item [c)]
The geodesic $c_n$ has no conjugate points in the interval
$[0,L_n)$;
\medskip
\item [d)]
For every $\ep>0$ there exists $T=T(\ep)>0$ such that for all
$n\in\N$ and $T<t<L_n-T$ we have $d(c_n(t),A)<\ep$.
\end{itemize}

\noindent Then $p$ and $q$ cannot be blocked away from each other
by a finite blocking set.
\begin{proof}
In the course of proof we will pass repeatedly from sequences of
geodesics to infinite subsequences. In order to avoid cumbersome notation,
we will use the following convention. Let $a_k$ be an infinite sequence.
We will denote its infinite subsequences by $a_k$ again, as opposed to,
say, $a_{n_k}$.

Assume that the claim fails. Then there is a point $z\in
M\setminus\{p,q\}$ and infinitely many $n\in\N$ such that
$z=c_n(t_n)$ for some $t_n\in(0,L_n)$. Passing to a subsequence,
we assume that for all $n\in\N$ there exist $t_n\in(0,L_n)$ such
that $z=c_n(t_n)$. Set
$$
\ep=\ep(p,q,z,A)=\frac12\min\left(\{d(p,A),d(q,A),d(z,A)\}\setminus\{0\}\right).
$$

Let $T=T(\ep)$ be as in condition d). Passing to a subsequence, if
need be, and using condition d), we obtain that either
$t_n\in(0,T)$ or $t_n\in(L_n-T,L_n)$ for all $n\in\N$.

Consider first the case $t_n\in(0,T)$. Again passing to a
subsequence, we assume that $\lim_{n\to\infty}t_n=\bar{t}$ exists,
and the geodesics $c_n$ converge to a geodesic $c:[0,\infty)\to
M$. Thus, we have
\begin{equation}   \label{points_eq}
c_n(t_n)=z=c(\bar{t})
\end{equation}
for all $n\in\N$. Since $c(0)=p$ and $z\ne p$, we conclude that
$\bar{t}>0$. By conditions a) and c), the geodesic $c$ has no
conjugate points. This contradicts equation~~\eqref{points_eq},
unless we have $c_n=c|_{[0,L_n]}$ for almost all $n\in\N$.

Suppose this is the case. Then, by condition a), there are
positive integers $n<m$ such that $T<L_n<L_m-T$ and
$q=c_n(L_n)=c(L_n)=c_m(L_n)$. Then, by condition b), $q\notin A$.
On the other hand, by condition d), we have $d(c_m(L_n),A)<\ep$.
This contradicts our choice of $\ep\le d(q,A)$. Thus,  we have
arrived at a contradiction.

\medskip

It remains to consider the case $t_n\in(L_n-T,L_n)$. Observe that
our setting is symmetric with respect to the interchange of $p$
and $q$ and simultaneous reversal of the time direction. This
symmetry flips the two cases at hand. Thus, the assumption
$t_n\in(L_n-T,L_n)$ leads to a contradiction as well.
\end{proof}
\end{prop}

\medskip

%\medskip

Our next result yields the existence of a nonempty open set of
insecure point pairs in a nonflat two-torus.

\begin{thm}    \label{no_block_thm}
Let $T^2$ be a nonflat riemannian two-torus. Let $Z\subset T^2$ be
an admissible cylinder, as introduced in section~~\ref{constr}.
Then every pair $(p,q)\in Z\times\overline{Z}$ is insecure.
\begin{proof}
We apply Proposition~~\ref{suffice_prop}  to our torus, with
$A=\bo Z$. Let $S$ be an admissible strip corresponding to $Z$;
let $P\in S,\,Q\in\overline{S}$ be such that $\pi(P)=p,\pi(Q)=q$.
Let $k\in\stab(S)$ be as in section~~\ref{constr}. For every
$n\in\N$ we choose a minimal geodesic
$\tc_n:[0,L_n]\to\overline{S}$ from $\tc_n(0)=P$ to
$\tc_n(L_n)=Q+nk$. Set $c_n=\pi\circ\tc_n$. By definition, this is
a sequence of geodesics $c_n:[0,L_n]\to T^2$ from $p$ to $q$.

Next, we check that the assumptions of
Proposition~~\ref{suffice_prop} are satisfied. Assumptions a) and
c) are satisfied by construction. Assumption b) is satisfied
because $p\in Z$ and the set $A=\bo Z$ is totally geodesic.
Assumption d) follows from Lemma~~\ref{key_lem} when $q\in Z$ and
from Lemma~~\ref{compl_lem} when $q\in\bo Z$.
\end{proof}
\end{thm}

\medskip

\begin{rem}     \label{ho_rem}
The preprint~~\cite{WKH} gives an argument to show the insecurity for pairs
$(p\in Z,q\in\bo Z)$. This observation is not contained in \cite{BaG08}
which is a preliminary version of the present paper.
\end{rem}
\medskip

\begin{thm}   \label{torus_thm}
A two-dimensional riemannian torus is secure if and only if it is
flat.
\begin{proof}
By \cite{Gut05} or \cite{GS06}, a flat torus is secure.
Theorem~\ref{no_block_thm} implies the existence of insecure point
pairs in any nonflat two-torus.
\end{proof}
\end{thm}

\medskip

We point out that the flat tori are distinguished amongst all riemannian two-tori
by the security of pairs $y=x$.
\begin{corol}    \label{flat_cor1}
A two-dimensional riemannian torus is flat if and only if all
pairs $(x,x)$ are secure.
\begin{proof}
It suffices to show that a nonflat riemannian two-torus contains
at least one point that cannot be blocked away from itself. Let
$Z$ be an admissible cylinder. Then any point $x\in Z$ cannot be
blocked away from itself.
\end{proof}
\end{corol}

\medskip

%Theorem~\ref{main_thm} implies that a nonflat riemannian two-torus always has insecure
%pairs of points. Are there any secure pairs in a nonflat two-torus? We think that there are none.
%Although our arguments do not prove this, they show that there is nonempty open set of insecure pairs.

%
%\begin{corol}    \label{flat_cor2}
%Let $(T^2,g)$ be a nonflat riemannian torus. Then the set of pairs of points $x,y$
%that cannot be blocked away from each other contains a nonempty open set.
%\begin{proof}
%Immediate from Proposition~\ref{no_block_prop}.
%\end{proof}
%\end{corol}
%

%
\section{Blocking and insecurity for tori of revolution}   \label{revolut}
By a {\em two-torus of revolution} we will mean the cartesian square  $T^2$
of the standard circle $T=\R/\Z$ with the riemannian metric
\begin{equation}   \label{tor_rev_eq}
ds^2=f^2(y)dx^2+dy^2.
\end{equation}
The somewhat more general $T$-invariant riemannian tori
\begin{equation}   \label{tor_inv_eq}
ds^2=f^2(y)dx^2+g^2(y)dy^2,
\end{equation}
reduce to the torus of revolution equation~~\eqref{tor_rev_eq} by a change of variables. 

\medskip

Thus, our torus of revolution is determined by the positive function $f:T\to\R_+$.
To emphasize the dependence on $f$ we will sometimes denote the
riemannian torus in equation~~\eqref{tor_rev_eq} by $T^2_f$.
We will assume that $f$ is strictly positive and infinitely differentiable.

Our analysis is based on the symmetries of $T^2_f$. The riemannian metric in equation~~\eqref{tor_rev_eq}
is invariant with respect to the action of $T$ on $T^2$ by $(x,y)\mapsto(x+a,y)$. 
Besides rotations $\rho_a(x)=x+a\mod1$, the isometry group of $T$ contains reflections.
With each $a\in T$ we associate the reflection $\si_a(x) = 2a-x\mod1$. Then $\si_a$ is
orientation reversing; it has two fixed points: $a,a+\frac{1}{2}$. 

\begin{prop}   \label{sub_main_thm}
Let $f:T\to\R_+$ be a smooth, positive function. Suppose that
i) it has a unique minimum, say at the point $a\in T$; ii) it is invariant under the reflection $\si_a$.

Then all pairs $\{(p,a),(q,a)\}\in T^2_f\times T^2_f$  are
secure. All other pairs of points in $T^2_f$ are insecure.
\begin{proof}
Let $F,G$ be transformations of $T$. We denote by
$F\times G$ the transformation of $T^2$ given by $(x,y)\mapsto(F(x),G(y))$. 
 
For arbitrary points $p,q\in T$ let $F_{p,q}:T\to T$ be the unique reflection
interchanging $p$ and $q$. If $r\in T$ is any of the two midpoints between $p$ and $q$,
then $F_{p,q}=\si_r$. We set
\begin{equation}    \label{symm_eq}
\ip_{p,q}=F_{p,q}\times\si_a. 
\end{equation}
Then $\ip_{p,q}:T^2\to T^2$ is an involution: $\ip_{p,q}^2=\id$; it has exactly $4$ fixed points. 

The coordinate vector fields $\bo/\bo x,\bo/\bo y$ on $T^2$ yield global coordinates on
the tangent bundle to $T^2_f$. Let $X$ be a tangent vector to $T^2$. The unique decomposition
$$
X=\xi\frac{\bo}{\bo x} + \eta\frac{\bo}{\bo y}
$$
defines the isomorphism $X\mapsto \left((x,y),(\xi,\eta)\right)$ of the tangent bundle  and $T^2\times\RR$.
If the base point of a tangent vector $X$ is clear from the context,
we will identify $X$ with $(\xi,\eta)\in\RR$.
Let $H:T^2\to T^2$ be a diffeomorphism. Our isomorphism identifies the differentials $D_{(x,y)}H$
with $2\times 2$ matrices. By  equation~~\eqref{symm_eq}, we have 
\begin{equation}    \label{diff_eq}
D_{(x,y)}\ip_{p,q}=-\id
\end{equation}
identically on $T^2$.

\medskip

From now on, we view $T^2$ as the riemannian torus of revolution $T^2_f$.
Let $y\in T$. By analogy with the round euclidean torus, we define the circles of latitude as
$L_y=T\times\{y\}\subset T^2_f$. 

\medskip

Let $X$ be the tangent vector with the base point $(x,y)$ and the representation 
$X=(\xi,\eta)$. The function on the tangent bundle given by
\begin{equation}    \label{clair_eq}
F(X)=f^2(y)\xi
\end{equation}
is the {\em Clairaut integral} for the geodesic flow on the tangent bundle of $T^2_f$  \cite{Kl}. 
Let $c:\R\to T_f^2$ be a geodesic. Set $c(t)=(x(t),y(t)),\dot{c}(t)=(\dot{x}(t),\dot{y}(t))$.
Let $t_1\ne t_2$ be such that $y(t_1)=y(t_2)$. Then, from equations~~\eqref{tor_rev_eq},~~\eqref{clair_eq}
\begin{equation}    \label{tangent_eq1}
\dot{x}(t_1)=\dot{x}(t_2),\ |\dot{y}(t_1)|=|\dot{y}(t_2)|.
\end{equation}

\medskip

Suppose now that $c(0)\in L_a$. Then, by equations~~\eqref{tor_rev_eq} and~~\eqref{clair_eq}
\begin{equation}    \label{minim_eq}
\min_{t\in\R} |\dot{y}(t)|=|\dot{y}(0)|.
\end{equation}
By equation~~\eqref{minim_eq}, the geodesics $c(t)=(x(t),y(t))$ in $T_f^2$ that intersect the circle
$L_a$ satisfy the following dichotomy: Either $\dot{y}(t)\equiv 0$ or $\dot{y}(t)$ does not change sign.
In the former case, the trace of $c$ is $L_a$. In the latter case, the geodesic $c$ winds around
$T_f^2$ crossing every circle of latitude transversally, as opposed to oscillating between two circles
of latitude. Thus, for the geodesics $c:\R\to T_f^2$ intersecting the circle
$L_a$, equation~~\eqref{tangent_eq1} becomes
\begin{equation}    \label{tangent_eq2}
\dot{x}(t_1)=\dot{x}(t_2),\ \dot{y}(t_1)=\dot{y}(t_2).
\end{equation}

\medskip

Let $\la$ be the free homotopy class of the circles of latitude. 
By condition i) on $f$, the closed curve $L_a\subset T^2_f$ is the unique minimal geodesic in $\la$. 
In view of Theorem~~\ref{no_block_thm}, it
suffices to show that any pair of points in $L_a$ is secure. To simplify notation, we will use the
shorthand $\{(p,a),(q,a)\}=\{p,q\}\in L_a\times L_a$.
For any torus of revolution $T^2_f$ the transformations $F_{p,q}\times\id:T^2_f\to T^2_f$ are isometries.
By assumption ii), the transformation $\id\times\si_a:T^2_f\to T^2_f$  is also an isometry of $T^2_f$.
Thus, the involutions $\ip_{p,q}:T^2_f\to T^2_f$ from  equation~~\eqref{symm_eq} are isometries.

\medskip

We denote by $G_{p,q}$ the set of geodesics $c:\R\to T^2_f$ such that for some
$t_p<t_q$ we have $c(t_p)=p,c(t_q)=q$. Restricting these geodesics to $[t_p,t_q]$ and taking their
traces, we obtain the set $\Ga_{p,q}$ of {\em geodesic segments connecting $p$ and $q$}. 
Note that $\Ga_{q,p}=\Ga_{p,q}$.
By equation~~\eqref{symm_eq}, $\ip_{p,q}:G_{p,q}\to G_{q,p}$, and hence $\ip_{p,q}(\Ga_{p,q})=\Ga_{p,q}$.
We will now show that actually $\ip_{p,q}$ maps every geodesic segment $c\in G_{p,q}$ to itself.

\medskip

Let $c\in G_{p,q}$. By equation~~\eqref{tangent_eq2}, $\dot{c}(t_p)=\dot{c}(t_q)$. In
view of equation~~\eqref{diff_eq}, every geodesic geodesic segment connecting $p$ and $q$
is invariant under $\ip_{p,q}$. Since $\ip_{p,q}$ interchanges $p$ and $q$, it fixes the 
midpoint of each $\ga\in\Ga_{p,q}$.

\medskip

Thus, the midpoint, say $m(\ga)$, of any $\ga\in\Ga_{p,q}$ is one of the four fixed points of $\ip_{p,q}$.
It may happen that $m(\ga)=p$. Then $\ga$ is a multiple of a prime geodesic segment
$\be\in\Ga_{p,q}$. Applying the same argument to $\be$, we conclude that  
the set of fixed points of $\ip_{p,q}$ minus 
the set $\{p,q\}$ is a blocking set for $\Ga_{p,q}$. 
\end{proof}
\end{prop}

\medskip

\begin{rem}      \label{midpoint_rem}
{\em 

A pair $\{p,q\}$ in a riemannian manifold is {\em midpoint secure} if
the set of midpoints of geodesics $\ga\in\Ga_{p,q}$ is finite. A manifold is
midpoint secure if every pair in $M$ is midpoint secure. Flat manifolds
are midpoint secure. Secure pairs in compact symmetric spaces are midpoint
secure \cite{GS06}. By our proof of Proposition~~\ref{sub_main_thm}, a pair
in a torus of revolution is either insecure or midpoint secure. These examples seem to suggest
that security is equivalent to the midpoint security. This is false, however. 
P. Herreros \cite{Her} gave an example of a riemannian two-sphere with a large family
of point pairs that are secure but not midpoint secure.  

}
\end{rem}

\medskip

\begin{corol}    \label{round_tor_cor}
Let $T^2=T^2(r,R)\subset\R^3$ be a round euclidean torus of revolution. See  
Example~~\ref{tor_rev_exa}. Let $E_{\inn}\subset T^2$ be the inner equator. Then
a point pair $\{p,q\}\in T^2\times T^2$ is secure if and only if $\{p,q\}\in E_{\inn}\times E_{\inn}$.
The pairs in $E_{\inn}$ are, in fact, midpoint secure.
\begin{proof}
We will refer to $T^2(r,R)\subset\R^3$ as the  round euclidean torus.
In the standard angular coordinates 
$0\le \al,\be \le 2\pi$ the round euclidean torus has the metric
\begin{equation}   \label{round_tor_rev_eq}
ds^2=(R-r\cos\be)^2\,d\al^2+r^2\,d\be^2.
\end{equation}
Since the round euclidean torus is a special case of the torus in equation~~\eqref{tor_inv_eq}, 
we reduce it to a
torus of revolution $T^2_f$ of equation~~\eqref{tor_rev_eq}. In this representation
we have $f(y)=R-r\cos(2\pi y)$,  up to a constant factor. This function has a unique minimum, $f(0)=R-r$,
and is invariant under the reflection $\si(y)=-y\mod1$. Hence, Proposition~~\ref{sub_main_thm} 
and Remark~~\ref{midpoint_rem} apply.
\end{proof}
\end{corol}

%\vspace{3mm}

%
\section{Insecurity for all point pairs in a two-dimensional torus}   \label{generic}
In the preceding section we saw that there are nonflat two-tori
with some secure point pairs. Now we will prove that such tori are
exceptional. We will exhibit a $C^2$-open and $C^{\infty}$-dense
subset $\G_{\tot}$ of the space of riemannian metrics
on the two-torus such that for metrics in $\G_{\tot}$ every pair
$(p,q)\in T^2\times T^2$ is insecure.

\medskip

We will recall a few well known facts about homotopically minimal
geodesics on an arbitrary riemannian two-torus $(T^2,g)$. First,
we explain what we mean by the {\em homological direction} of a
homotopically minimal geodesic $c:\R\to T^2$. If $\tc:\R\to\RR$ is
a lift of $c$, then there is a strip bounded by two parallel lines
in $\RR$ that contains $\tc$. The minimality of $c$ implies that
$\tc$ is proper and separates the boundary components of the
strip. This was first proven by G.A. Hedlund \cite{He}. See also
\cite{Ba1}, Theorem 6.5. Moreover, $\tc$ determines an orientation
for the parallel lines in the strip. This oriented direction can
be viewed as an element in the quotient of $\RR\setminus\{0\}$ by
the multiplicative action of $\R_+$. Representing $T^2$ as
$\RR/\ZZ$, as usual, we identify $\ZZ\subset\RR$ with
$H_1(T^2,\Z)$ and $\RR$ with $H_1(T^2,\R)$.

The homological direction $[h]\in
(H_1(T^2,\R)\setminus\{0\})/\R_+$ of the homotopically minimal
geodesic $c$ is the oriented direction of a strip that contains a
lift of $c$. Let $c$ be a periodic geodesic of minimal length in
its free homotopy class. Then $c$ is homotopically minimal, by
Theorem~~\ref{minim_length_thm}. Let $h=h(c)$ be its homology
class. Then the homological direction of $c$ is the image of
$h(c)\in H_1(T^2,\Z)$ in the quotient
$(H_1(T^2,\R)\setminus\{0\})/\R_+$.

Alternatively, one can define the homological direction of a
homotopically minimal geodesic as follows. 
Let $c:\R\to T^2$ be
one. Let $s_i,t_i\in\R,\,i\in\N,$ be two infinite sequences such
that $s_i<t_i$ and $t_i-s_i\to\infty$. Let $\ga_i,i\in\N,$ be
arbitrary curves from $c(t_i)$ to $c(s_i)$ such that their lengths
are uniformly bounded. We will denote by $\al*\be$ the concatenation
of curves $\al$ and $\be$ whenever it is defined. Let $h_i\in H_1(T^2,\R)$ be the homology
class of the closed curve $c|_{[s_i,t_i]}*\ga_i$. Then, for
all but a finite number of indices, we have $h_i\ne 0$, defining
the direction $[h_i]\in (H_1(T^2,\R)\setminus\{0\})/\R_+$.
Moreover, as $i\to\infty$, the sequence $[h_i]$ converges; the limit
$[h]=\lim_{i\to\infty}[h_i]$ is the homological direction of $c$.
See \cite{Sch} for a similar concept.

\medskip

Let $S$ be an oriented surface. Let $v_1,v_2$ be linearly
independent tangent vectors with the same base point in $S$. We
set $\ep(v_1,v_2)=1$ if the pair $v_1,v_2$ is positively oriented,
and $\ep(v_1,v_2)=-1$ otherwise. Let $[h_1],[h_2]\in
(H_1(T^2,\R)\setminus\{0\})/\R_+$. Let $h_i\in[h_i]$ for $i=1,2$.
Then the sign of the homological intersection $h_1\cdot h_2$ does
not depend on the choices of $h_i\in[h_i]$. Hence, we will use the
notation $\sgn\,(h_1\cdot h_2)\in\{0,\pm 1\}$.

We will need the following property of homotopically minimal
geodesics on a two-torus.
\begin{lem}      \label{intersect_lem}
Let $c_1:\R\to T^2,\,c_2:\R\to T^2$ be homotopically minimal
geodesics. Let $[h_1],[h_2]\in (H_1(T^2,\R)\setminus\{0\})/\R_+$
be their homological directions.

\noindent Suppose that $\sgn\,(h_1\cdot h_2)\ne 0$. Then for any
$t_1,t_2\in\R$ such that $c_1(t_1)=c_2(t_2)$, we have
$$
\ep(\dot{c}_1(t_1),\dot{c}_2(t_2))=\sgn\,(h_1\cdot h_2).
$$
\begin{proof}
Let $\tc_i:\R\to\RR,i=1,2,$ be the lifts of $c_i$ such that
$\tc_1(t_1)=\tc_2(t_2)$. We will prove that
$\ep(\dot{\tc}_1(t_1),\dot{\tc}_2(t_2))=\sgn\,(h_1\cdot h_2)$.

We identify $H_1(T^2,\R)$ with $\RR$, as always, and set
$l_i=\{th_i:t\in\R\}$ where $i=1,2$. Note that $l_1,l_2$ are
oriented straight lines in $\RR$. For $(s,t)\in\R\times[0,1]$ and
$i=1,2$,  set $H_i(s,t)=(1-t)\tc_i(s)+tsh_i$. Then
$H_i:\R\times[0,1]\to\RR$ is a homotopy between $\tc_i$ and the
line $l_i$. Set $K=\{(s,t)\in\R\times[0,1]: H_1(s,t)\in l_2\}$. By
remarks in the beginning of this section,  $K$ is a compact subset
of $\R\times[0,1]$. Standard arguments from differential
topology\footnote{See, e. g., \cite{Hi}, Chapter 5.} show that the
oriented intersection numbers
$\#(\tc_1,\tc_2),\#(\tc_1,l_2),\#(l_1,\tc_2)$ are defined and
satisfy
\begin{equation}    \label{intersect_eq1}
\#(\tc_1,\tc_2)=\#(\tc_1,l_2)=\#(l_1,l_2).
\end{equation}
We obviously have
\begin{equation}    \label{intersect_eq2}
\#(l_1,l_2)=\sgn\,(h_1\cdot h_2).
\end{equation}

\medskip

Since $\tc_1,\tc_2$ are minimal geodesics, they intersect exactly
once, and transversely, at the point $\tc_1(t_1)=\tc_2(t_2)$.
Hence
\begin{equation}    \label{intersect_eq3}
\#(\tc_1,\tc_2)=\ep(\dot{\tc}_1(t_1),\dot{\tc}_2(t_2)).
\end{equation}
Equations~~\eqref{intersect_eq1} -~~\eqref{intersect_eq3} imply
our claim.
\end{proof}
\end{lem}

\medskip

\begin{corol}      \label{intersect_cor}
Let $h_1,h_2$ be prime elements in the lattice of integer classes
in $H_1(T^2,\R)$. Let $c_1,c_2:\R\to T^2$ be minimal periodic
geodesics in the classes $h_1,h_2$ respectively. Suppose that
$h_1\cdot h_2\ne 0$. Then
$$
\card(c_1(\R)\cap c_2(\R))=|h_1\cdot h_2|.
$$
\begin{proof}
Denote by $L_i$ the minimal period of $c_i$. Since $c_1$ and $c_2$
intersect transversely, we have
$$
h_1\cdot h_2=\sum_{(t_1,t_2)}\ep(\dot{c}_1(t_1),\dot{c}_2(t_2))
$$
where the sum is over all pairs $(t_1,t_2)\in[0,L_1)\times[0,L_2)$
such that $c_1(t_1)=c_2(t_2)$. By Lemma~~\ref{intersect_lem}
\begin{equation}    \label{intersect_eq4}
|h_1\cdot
h_2|=\card\{(t_1,t_2)\in[0,L_1)\times[0,L_2):\,c_1(t_1)=c_2(t_2)\}.
\end{equation}

Since $h_1,h_2\in H_1(T^2,\R)$ are prime, the curves
$c_1:[0,L_1)\to T^2$ and $c_2:[0,L_2)\to T^2$ are injective.
Hence, for every $p\in c_1(\R)\cap c_2(\R)$ there is exactly one
pair $(t_1,t_2)\in[0,L_1)\times[0,L_2)$ such that
$p=c_1(t_1)=c_2(t_2)$. Now our claim follows directly from
equation~~\eqref{intersect_eq4}.
\end{proof}
\end{corol}

\medskip

We will investigate riemannian metrics on the two-torus. We denote
by $\G_1,\G_2$ the sets of riemannian metrics on $T^2$ satisfying
the conditions $G_1,G_2$ respectively:

\medskip

\begin{itemize}
\item [$(G_1)$] There exists a homology class $h\in H_1(T^2,\R)\setminus\{0\}$ such that
$T^2$ does not admit a foliation by homotopically minimal geodesics with the
homological direction $[h]$;\\
\medskip
\item [$(G_2)$] There exist integral homology classes $h_1,h_2\in
H_1(T^2,\R)$ such that $h_1\cdot h_2=1$ and for each $i=1,2$ there is a unique
minimal periodic geodesic in the class  $h_i$.

\end{itemize}

\medskip

We will now show that for $g\in\G_1\cap\G_2$ any point pair in
$(T^2,g)$ is insecure.

\begin{thm}     \label{total_insec_thm}
Let $(T^2,g)$ be a riemannian torus. If the metric $g$ satisfies
conditions $(G_1)$ and $(G_2)$, then every pair of points in
$(T^2,g)$ is insecure.
\begin{proof}
Suppose the claim is false, and let $(p,q)$ be a secure pair of
points in our torus. Let $h_1,h_2\in H_1(T^2,\R)$ be as in
condition $(G_2)$; we denote by $c_1\in h_1,c_2\in h_2$ the unique
minimal periodic geodesics in these classes. Then
$Z_i=T^2\setminus c_i(\R)$ are admissible cylinders. By
Theorem~~\ref{no_block_thm}, both points $p,q$ belong to the sets
$c_1(\R),c_2(\R)$. Thus, $\{p,q\}\subset c_1(\R)\cap c_2(\R)$.

By Corollary~~\ref{intersect_cor}, the geodesics $c_1,c_2$
intersect at a unique point. Thus $p=q$, and $\{p\}=c_1(\R)\cap
c_2(\R)$. We will now study the geodesics in $(T^2,g)$ passing
through $p$.

\medskip

In view of Theorem~~\ref{no_block_thm}, we know that for any free
homotopy class of prime periodic geodesics, there exists a minimal
geodesic in this class passing through $p$. Taking limits of these
geodesics, we conclude that for every $h\in
H_1(T^2,\R)\setminus\{0\}$ there exists a homotopically minimal
geodesic with the homological direction $[h]$ passing through $p$.
By Theorem~~6.7 and Theorem~~6.9 in \cite{Ba1}, this geodesic is
uniquely determined by $[h]$. We denote it by $c([h])$. Let
$UT_p^2\subset T_p^2$ be the set of unit tangent vectors at $p$.
Let $\Phi([h])\in UT_p^2$ be the tangent vector at $p$ to the
geodesic $c([h])$. This defines a map
$$
\Phi:(H_1(T^2,\R)\setminus\{0\})/\R_+\to UT_p^2.
$$

\medskip

By Theorem~~VII in \cite{He}, or by Remark $(2)$ on page $32$ in
\cite{Ba1}, there is a uniform bound on the widths of the strips
in $\RR$ containing lifts of homotopically minimal geodesics. It
follows that the map $\Phi$ is continuous. Since it is also
injective, we conclude that $\Phi$ is a homeomorphism. Therefore,
every geodesic in $(T^2,g)$ passing through $p$ is homotopically
minimal.

\medskip

We have shown that $p$ is a {\em pole} for the riemannian torus
$(T^2,g)$. This means that any geodesic passing through $p$ has no
points conjugate to $p$. Now we use the following result of
\cite{Ba4}: If $(T^2,g)$ has a pole, then for every $[h]\in
(H_1(T^2,\R)\setminus\{0\})/\R_+$ the torus admits a foliation  by
minimal geodesics with the homological direction $[h]$. Thus, the
riemannian torus $(T^2,g)$ does not satisfies condition $(G_1)$.
We have arrived at a contradiction.
\end{proof}
\end{thm}

\medskip

We have shown that for riemannian tori $(T^2,g)$, where
$g\in\G_1\cap\G_2$, every pair of points is insecure. Recall that
a closed geodesic $\ga$ is {\em nondegenerate} if the linearized
Poincar\'e map for $\ga$ does not have eigenvalue $1$.

Let $\G_3\subset\G_2$ be the set of tori $(T^2,g)$ satisfying the
following strengthening of condition $(G_2)$:
\begin{itemize}
\item [$(G_3)$] There exist integral homology classes $h_1,h_2\in
H_1(T^2,\R)$ such that $h_1\cdot h_2=1$ and there is a unique
minimal periodic geodesic in the class  $h_i,\,i=1,2$. Moreover,
these geodesics are nondegenerate.
\end{itemize}

\medskip

%Our next proposition states that $\tG$ is a generic set of
%metrics.
%
\begin{prop}     \label{generic_prop}
The set $\G_1\cap\G_3$ of riemannian metrics on the two-torus
satisfying conditions $(G_1)$ and $(G_3)$ is $C^2$-open and
$C^{\infty}$-dense in the space of riemannian metrics.
\end{prop}
\medskip

Our proof of Proposition~~\ref{generic_prop} is based on the
following two lemmas.
\begin{lem}       \label{condit1_lem}
The set $\G_1$ of metrics on the two-torus satisfying condition $(G_1)$
is $C^2$-open and $C^{\infty}$-dense.
\begin{proof}
Let $h\in H_1(T^2,\R)\setminus\{0\}$ be an integral class. Denote by $\G_1(h)$ the set
of metrics that do not admit a foliation by homotopically minimal
geodesics with homological direction $[h]$. We will prove first that $\G_1(h)$
is $C^{\infty}$-dense. Since $\G_1(h)\subset\G_1$, this will prove that 
$\G_1$ is $C^{\infty}$-dense.

Let $g$ be a metric in the complement  $\G_1(h)^c$ of $\G_1(h)$, i. e., 
there exists a foliation $\FF$ of $T^2$ by homotopically minimal
geodesics (with respect to $g$) of homological direction $[h]$.
Then $\FF$ has the following structure: i) The torus contains a nonempty 
closed subset  foliated by periodic geodesics in  $\FF$;
ii) The complement to this set is either empty or it is a disjoint
union of admissible cylinders foliated by geodesics in $\FF$
that are heteroclinic to the boundary geodesics of the cylinder.
See, e. g., \cite{Ba1}, Theorem~~6.6 and Theorem~~6.7.

We will now show that if not all of $T^2$ is foliated by periodic geodesics
in  $\FF$, then $g$ can be  $C^{\infty}$-approximated by a sequence in $\G_1(h)$.
Let $Z\subset T^2$ be an admissible
cylinder foliated by  heteroclinic geodesics in $\FF$. We choose two minimal
heteroclinic geodesics $c_-:\R\to Z$ and $c_+:\R\to Z$, one for each of the two possibilities;
see Lemma~~\ref{component_lem}. Then $U=Z\setminus(c_-(\R)\cup c_+(\R))$ is a 
nonempty open subset of $Z$. We choose $p\in U$ and a sequence $g_k$ of riemannian metrics
such that a) $g_k\ge g$, b) $g_k|_p>g|_p$, and c) $g_k|_{T^2\setminus U}=g|_{T^2\setminus U}$.
It follows that each of the tori $(T^2,g_k)$ has neither a heteroclinic
nor a periodic homotopically minimal geodesic with homological direction
$[h]$ passing through $p$. Hence $g_k\in\G_1(h)$ for all $k\in\N$. This proves that
$\G_1(h)$ is $C^{\infty}$-dense in the set of all metrics on $T^2$ that do
not admit a foliation by periodic minimal geodesics of direction $[h]$.

It remains to prove that every metric that admits such a foliation can be
$C^{\infty}$-approximated by metrics that do not admit such a foliation.
This can be shown by a perturbation argument, which is
similar, but simpler, than the argument in the first part of the proof.

\medskip

To prove that $\G_1$ is open in the $C^2$-topology, it suffices to
show that $\G_1(h)^c$ is $C^2$-closed for every $h\in H_1(T^2,\R)\setminus\{0\}$. 
This fact is closely related to the convergence
properties for invariant circles of monotone twist maps. 
For the benefit of the reader, we will outline a proof.

Let $g_i\in\G_1(h)^c$ be an infinite sequence such that $g_i\to g$
in the $C^2$-topology. Let $\FF_i,\,i=1,2,\dots,$ be the
corresponding foliations of $T^2$. Denote by $U_i$ the vector
field formed by the unit tangent vectors on $(T^2,g_i)$ tangent to
the foliation $\FF_i$. Then the sequence $U_i$ is equicontinuous, since
the vector fields $U_i$ are uniformly lipschitz. See \cite{So} or
\cite{Mat}. By the Arzela-Ascoli theorem, there is a subsequence
of vector fields $U_i$ that converges to a $g$-unit tangent vector
field $U$. The integral curves of $U$ form a foliation $\FF$ of
$(T^2,g)$ by homotopically $g$-minimal geodesics. The  widths of the parallel
strips containing the lifts of the geodesics in the foliations
$\FF_i$ to $\RR$ can be chosen uniformly bounded. See, e. g., the
proof of Remark $(2)$ on p. $32$ in
\cite{Ba1}.\footnote{Alternatively, we can use the {\em order of
the $\ZZ$-translates} of such lifts to determine the direction.
See, e. g., \cite{We} or \cite{Ba3}.} It follows that the
geodesics in $\FF$ have homological direction $[h]$. Hence,
$g\in\G_1(h)^c$.
\end{proof}
\end{lem}

\medskip

\begin{lem}       \label{condit23_lem}
The set $\G_3$ of riemannian metrics on the two-torus satisfying
condition $(G_3)$ is $C^2$-open and $C^{\infty}$-dense.
\begin{proof}
First, we will show that $\G_3$ is $C^2$-open. Suppose that $g$
satisfies condition $(G_3)$ for integer classes $h_1,h_2\in
H_1(T^2,\R)$ with $h_1\cdot h_2=1$. Let $c_1,c_2$ be the unique
$g$-minimal periodic geodesics in the classes $h_1,h_2$
respectively. Let $g_k,\,k\ge 1,$ be a sequence of riemannian
metrics on $T^2$ converging to $g$ in the $C^2$-topology. We will
show that all but a finite number of the metrics $g_k$ satisfy
condition $(G_3)$ with the classes $h_1,h_2$.

\medskip

Suppose first that an infinite subset of these metrics has several
minimal periodic geodesics in one of the homology classes, say
$h_1$. Relabelling the indices, we obtain a sequence of riemannian
metrics $g_k$ converging to $g$ and such that there exist distinct
$g_k$-minimal periodic geodesics $c_k'$ and $c_k''$ in the class
$h_1$. They are distinct in the sense that $c_k'(\R)\ne
c_k''(\R)$. With an appropriate reparameterization, both sequences
of geodesics converge to the $g$-minimal geodesic $c_1$. By the
implicit function theorem, the linearized Poincar\'e map for $c$
has an eigenvalue $1$, contrary to our assumption.

\medskip

Deleting a finite number of indices, we assume without loss of
generality that for all $k\in\N$ the metrics $g_k$ have unique
minimal periodic geodesics $c_k^{(1)},c_k^{(2)}$ in the classes
$h_1,h_2$ respectively. We claim that for all but a finite number
of indices the geodesics $c_k^{(1)},c_k^{(2)}$  are nondegenerate.
We may assume that we have the convergence $c_k^{(1)}\to
c_1,c_k^{(2)}\to c_2$. Since the geodesic flows of $(T^2,g_i)$
converge in the $C^1$-topology to the geodesic flow of $(T^2,g)$,
and the limit geodesics are nondegenerate, we obtain the claim.

\medskip

We have shown that $\G_3$ is $C^2$-open. We will now prove that
$\G_3$ is $C^{\infty}$-dense. Let $g$ be any $C^{\infty}$ metric
on the torus. Let $h_1,h_2$ be any integer homology classes such
that $h_1\cdot h_2=1$. Let $c_1\in h_1,c_2\in h_2$ be $g$-minimal
periodic geodesics.

Let $f$ be a nonnegative $C^{\infty}$ function on $T^2$ satisfying
the following conditions:
\begin{itemize}
\item[$a)$] We have $f^{-1}(0)=c_1(\R)\cup c_2(\R)$;

\medskip

\item[$b)$] For $i=1,2$ there exists a unit tangent vector $v_i$ for the torus $(T^2,g)$ with
base point in $c_i(\R)$ which is normal to $c_i(\R)$  and such
that $(f\circ\ga)''(0) > 0$ for every curve $\ga$ in $T^2$ with $\dot{\ga}(0)=v_1$ or $\dot{\ga}(0)=v_2$. 

\end{itemize}

\medskip

\noindent See figure~~\ref{condit23_fig}. Let $\la>0$. We define
the riemannian torus $(T^2,g_{\la})$ by $g_{\la}=(1+\la f)g$. By
condition $a)$, the curves $c_i:\R\to T^2$ are the unique
$g_{\la}$-minimal periodic geodesics in the classes $h_i, i=1,2$.
We claim that for all $\la>0$ the $g_{\la}$-geodesics $c_1,c_2$
are nondegenerate.

\begin{figure}[htbp]
\begin{center}
\input{condit23.pstex_t}
\caption{Illustration to the proof of Lemma~~\ref{condit23_lem}.}
\label{condit23_fig}
\end{center}
\end{figure}

\medskip

Assume the opposite, i. e., that for some $\la>0$ one of them, say
$c_1$, is a degenerate $g_{\la}$-geodesic. Then there exists a
nontrivial, periodic, normal $g_{\la}$-Jacobi field $Y$ along
$c_1$. Let $\exp_{\la}$ be the exponential map of $g_{\la}$. For
$\tau\in\R$ denote by $\al_{\tau}(t)=\exp_{\la}(\tau Y(t))$ the
normal variation of $c_1$ induced by $Y$. Let $L_{\la}$ denote the
length with respect to $g_{\la}$; let $L$ be the length with
respect to $g$. The formula for the second variation of
arc-length\footnote{See, for instance, \cite{GKM}, page 122, for
the second variation of arclength formula.} implies
\begin{equation}    \label{deriv_eq}
\frac{d^2}{d\tau^2}L_{\la}(\al_{\tau})|_{\tau=0}=0.
\end{equation}

\medskip

On the other hand, since the Jacobi field $Y$ has only isolated
zeros, condition b) implies that for a nonempty open interval of
values of $t$ we have
\begin{equation}    \label{deriv1_eq}
\frac{d^2}{d\tau^2}f(\exp_{\la}(\tau Y(t)))|_{\tau=0}>0.
\end{equation}
Since $c_1$ is a $g$-minimal periodic geodesic, we have
\begin{equation}    \label{deriv2_eq}
\frac{d^2}{d\tau^2}L(\al_{\tau})|_{\tau=0}\ge 0.
\end{equation}

Differentiating the expression
\begin{equation}    \label{length_eq}
L_{\la}(\al_{\tau})=\int_0^{L(c_1)}\left[\left(1+\la f(\exp_{\la}(\tau
Y(t)))\right)g\left(\dot\al_{\tau}(t),\dot\al_{\tau}(t)\right)\right]^{\frac12}dt
\end{equation}
and using equations~~\eqref{deriv1_eq} and~~\eqref{deriv2_eq}, we
obtain
$$
\frac{d^2}{d\tau^2}L_{\la}(\al_{\tau})|_{\tau=0}> 0.
$$

\medskip

This contradicts equation~~\eqref{deriv_eq}, and hence proves our
claim. Thus, $g_{\la}\in\G_3$ for all $\la>0$. By construction,
$(T^2,g_{\la})$ converges to $(T^2,g)$ as $\la\to 0$. Since $g$ is an
arbitrary metric, this completes the proof.
\end{proof}
\end{lem}

\medskip

\noindent{\em Proof of Proposition~~\ref{generic_prop}}. We recall
that $\G_1$ (resp. $\G_3$) is the set of riemannian tori $(T^2,g)$
satisfying condition $(G_1)$ (resp. $(G_3)$). By
Lemmas~~\ref{condit1_lem} and~~\ref{condit23_lem}, $\G_1$ and
$\G_3$ are $C^2$-open, hence the set $\G_1\cap\G_3$ is $C^2$-open
as well. 

By Lemmas~~\ref{condit1_lem} and~~\ref{condit23_lem}, the sets
$\G_1$ and $\G_3$ are $C^{\infty}$-dense. Since the $C^{\infty}$-topology is stronger than the
$C^2$-topology, $\G_1$ and $\G_3$ are also $C^{\infty}$-open.
The intersection of two
open and dense sets is open and dense. Thus, the set
$\G_1\cap\G_3$ is $C^{\infty}$-dense. \qed

\medskip

Theorem~~\ref{total_insec_thm} and Proposition~~\ref{generic_prop}
immediately imply the following.
\begin{corol}    \label{open_dense_cor}
The set of riemannian  tori $(T^2,g)$ contains a $C^2$-open and
$C^{\infty}$-dense subset of totally insecure tori.
\begin{proof}
Set $\G_{\tot}=\G_1\cap\G_3$. By Proposition~~\ref{generic_prop},
$\G_{\tot}$ is a $C^2$-open and $C^{\infty}$-dense subset of the
set of riemannian tori. By Theorem~~\ref{total_insec_thm}, for
$g\in\G_{\tot}$, every pair of points in $(T^2,g)$ is insecure.
\end{proof}
\end{corol}
%

%\vspace{5mm}

%
\section{Total insecurity for surfaces of genus greater than one}   \label{higher}

\medskip

In this section we will prove that on every closed riemannian
surface $M$ of genus greater than one every point pair is
insecure. Replacing $M$ by a two-sheeted covering space, if
necessary, we can assume without loss of generality that $M$ is
oriented. The basic idea is simple and has already been used to
prove the total insecurity of admissible cylinders. See
Theorem~~\ref{no_block_thm}.

\medskip

Besides the given metric $g$ on $M$, we consider a ``background metric''
$g_0$ of curvature  $-1$ on $M$. Then the universal riemannian covering 
of $(M,g_0)$ can be identified with the open unit disc 
$D\subseteq\mathbb C$ with the standard hyperbolic metric; thus, we have 
a covering $\pi:D\to M$ such that $\tg_0=\pi^*g_0$ is the  
standard hyperbolic metric. We will say that the
geodesics, distances, etc, with respect to $g_0$ or $\tg_0$ are
the hyperbolic geodesics, distances, etc. Note that $\tg_0$ and the lifted metric
$\tg=\pi^*g$ are lipshitz equivalent.

\medskip

The proof of Theorem~~\ref{high_genus_thm}, which is the main
result in this section, crucially uses fundamental facts due to
M.~Morse~~\cite{Mo}. For convenience of the reader, we give their
complete statements below.

\medskip

Theorem~~\ref{bound_dist_thm} invokes the standard notion of the
Hausdorff distance between subsets of a metric space. We do not
specify the metric in the statement of the theorem because
lipschitz equivalent metrics yield equivalent Hausdorff distances.
By a geodesic segment we will mean the trace of a geodesic, when
the endpoints are specified. 

\medskip

\begin{thm}     \label{bound_dist_thm}
(\cite{Mo}, Lemma 8). There is a positive constant $R$ such that
the Hausdorff distance between any $\tg$-minimal geodesic segment
and the hyperbolic geodesic segment with the same endpoints is
less than $R$.
\end{thm}

\medskip

In view of Theorem~~\ref{bound_dist_thm}, every $\tg$-minimal
geodesic $c:{\mathbb R}\rightarrow D$ is a uniformly bounded
distance away from some hyperbolic geodesic. Therefore, there are
distinct points $c(-\infty),c(\infty)\in
\partial D$ such that we have
$\lim\limits_{t\rightarrow \pm \infty}
c(t)=c(\pm\infty)$.\footnote{The limits are with respect to the
canonical topology on ${\mathbb C}$.} The statement below directly
follows from the fact that complete minimal geodesics intersect at
most once.

\begin{thm}     \label{no_intersect_thm}
(\cite{Mo}, Theorem 4). Let $c_1:{\mathbb R}\rightarrow D$ and
$c_2:{\mathbb R}\rightarrow D$ be $\tg$-minimal geodesics. If
$c_2(-\infty)$ and $c_2(\infty)$ belong to the same component of
$\partial D\setminus\{c_1(-\infty), c_1(\infty)\}$ then the
geodesics $c_1$ and $c_2$ do not intersect.
\end{thm}

\medskip

Our next  statement directly follows from
Theorem~~\ref{no_intersect_thm}.

\begin{thm}              \label{one_to_one_thm}
Let $c:\R/L\Z\rightarrow M$ be a periodic geodesic of minimal
$g$-length in its free homotopy class.  If this class contains a
simple closed curve, then the mapping $c:[0, L)\to M$ is
one-to-one.
\end{thm}
\begin{thm}              \label{two_minim_thm}
(\cite{Mo}, Theorem 11). Let $z, w\in\partial D$ be distinct
points.

\medskip

\noindent 1. Suppose that there is more than one $\tg$-minimal
geodesic $\tc:\R\to D$ with $\tc(-\infty)=z$ and $\tc(\infty)=w$.
Then there exist two $\tg$-minimal geodesics $\tc_{\pm}:\R \to D$
satisfying
$$
\tc_{\pm}(-\infty)=z, \tc_{\pm}(\infty)=w
$$
such that their traces bound a closed strip $S\subset D$, $S
\simeq\R\times[0,1]$, containing every $\tg$-minimal geodesic with
endpoints $z$ and $w$.

\noindent 2. Suppose, in addition, that $z$ and $w$ are the
endpoints of a lift of a closed curve in $M$. Let $\al$ be the
free homotopy class of this curve. Then $\tc_-$ and $\tc_+$
project to periodic geodesics of minimal $g$-length in the class
$\alpha$.
\end{thm}

\medskip

Let $\al$ be a nontrivial free homotopy class of closed curves in
$M$ that contains a simple closed curve. We will construct a
$\tg$-convex set $E\subseteq D$. In our proof of
Theorem~~\ref{high_genus_thm} the set $E$ will be a counterpart of
the strip $S\subseteq\RR$ in the proof of
Theorem~~\ref{no_block_thm} on the insecurity of nonflat tori.
There are two cases to consider.

\begin{figure}[htbp]
\begin{center}
\input{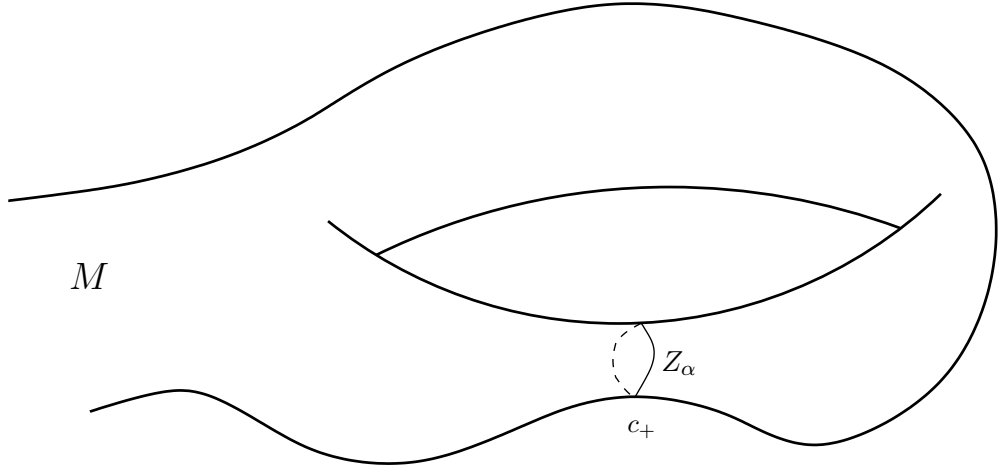}
\caption{Case (i) : there is a unique closed geodesic of minimal
length in the homotopy class $\alpha$.} \label{fig_high1}
\end{center}
\end{figure}

\medskip

\noindent Case (i): There is a unique periodic geodesic, say
$c_+$, of minimal $g$-length in the homotopy class $\al$. We
denote by $Z_{\al}\subset M$ the trace of $c_+$.

\medskip

\begin{figure}[htbp]
\begin{center}
\input{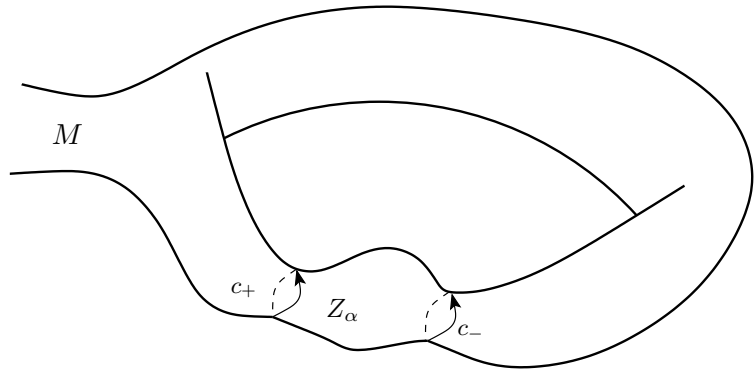}
\caption{Case (ii) : there are at least two closed geodesics of
minimal length in the homotopy class $\alpha$.} \label{fig_high2}
\end{center}
\end{figure}

\noindent Case (ii): There are at least two periodic geodesics of
minimal $g$-length in the class $\al$. It follows from
Theorem~~\ref{one_to_one_thm} and Theorem~~\ref{two_minim_thm}
that there exist two such geodesics bounding a closed cylinder
$Z_{\al}\subset M$ that contains every periodic geodesic of
minimal $g$-length in the class $\al$. We denote these geodesics
by $c_+$ and $c_-$, so that $Z_{\al}$ is to the left of $c_-$.
Figure~~\ref{fig_high1} and figure~~\ref{fig_high2} illustrate the
cases (i) and (ii) respectively. In particular, $Z_{\al}\simeq
S^1$ in case (i) and $Z_{\al}\simeq S^1\times[0,1]$ in case (ii).

\medskip

\begin{figure}[htbp]
\begin{center}
\input{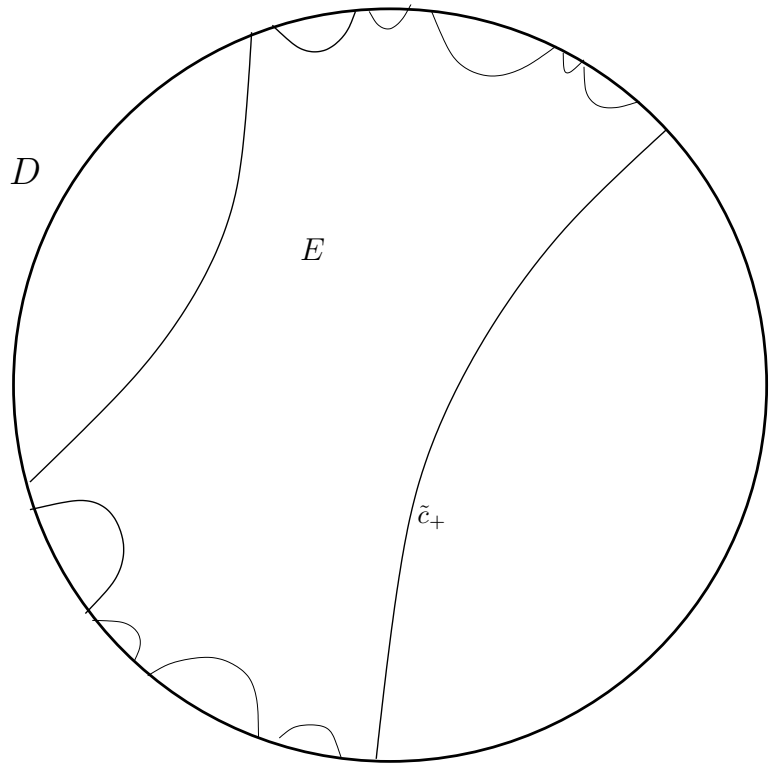}
\caption{The region $E$ in the case when there is a unique closed
geodesic of minimal length in the homotopy class $\alpha$.}
\label{fig_high3}
\end{center}
\end{figure}

\begin{figure}[htbp]
\begin{center}
\input{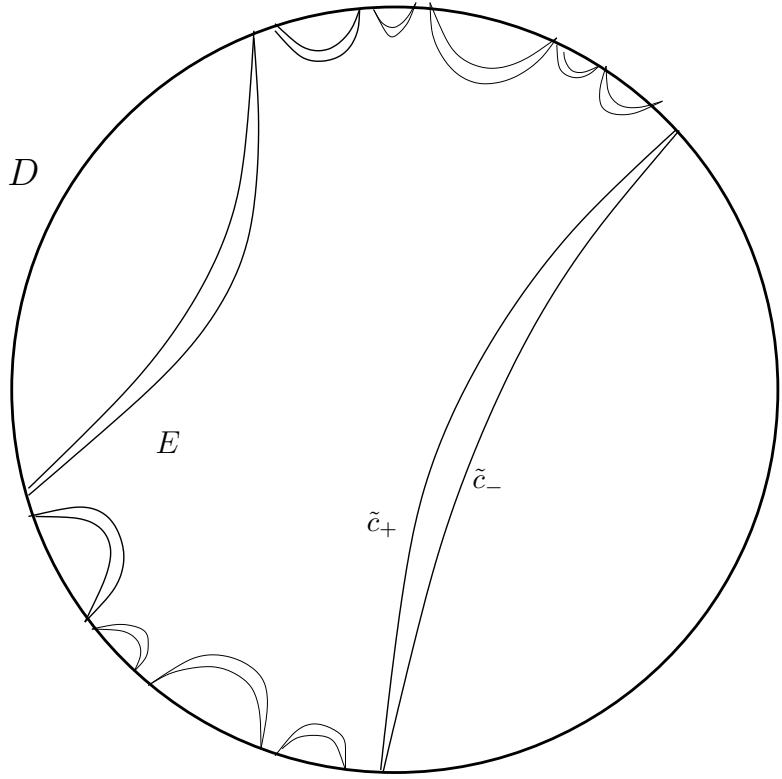}
\caption{The region $E$ in the case when there are at least two
closed geodesics of minimal length in the homotopy class
$\alpha$.} \label{fig_high4}
\end{center}
\end{figure}

We choose a lift, say $\tc_+$, of $c_+$ to $D$. Now we define $E$
to be the connected component of $D\setminus\pi^{-1}(Z_{\al})$
whose boundary contains $\tc_+$ and which lies to the left of
$\tc_+$. Note that the latter condition is automatically satisfied
in case (ii). See figure~~\ref{fig_high3} and
figure~~\ref{fig_high4} for illustration.

\vspace{2mm}

%

%We will now establish several properties of the set $E\subset D$.

%
\begin{lem}              \label{convex_lem}
The set $E$ has the following properties:
\begin{itemize}
\item[a)]
If $x,y\in E$ then every $\tg$-minimal geodesic from $x$ to $y$ is
contained in $E$.
\item[b)]
Let $\tal$ denote the covering transformation for
$\pi:D\rightarrow M$ mapping $\tc_+(0)$ to $\tc_+(L(c_+))$. Then
$\tal(E)=E$.
\item[c)]
Suppose, in addition, that $\alpha$ is homologically nontrivial.
Then $\pi(E)=M\setminus Z_{\al}$.
\end{itemize}
\begin{proof}
a) Assume the opposite. Let $c:[0, L]\to D$ be a $\tg$-minimal
geodesic  from $x$ to $y$ intersecting $\bo E$. Let $t_0\in(0,L)$
be the first time that $c$ intersects $\bo E$. Then there exists a
lift $\tc:\R\to D$ of $c_+$ or $c_-$ such that $c(t_0)\in\tc(\R)$
and $c([0, t_0))\subset E$.

Since $\tc(\R)$ disconnects $D$ and since the geodesics $c$ and
$\tc$ intersect transversely, there exists $t_1\in(t_0, L)$ such
that $c(t_1)\in\tc(\R)$. By Theorem~~\ref{minim_length_thm}, the
geodesic $\tc$ is $\tg$-minimal. Hence, the $\tg$-minimal geodesic
$c$ intersects $\tc$ at most once.

\vspace{2mm}

\noindent b) Since $\tal$ is a deck transformation, it preserves
any set of the form $\pi^{-1}(X),\,X\subset M$. In particular,
$\tal\left(\pi^{-1}(Z_{\al})\right) = \pi^{-1}(Z_{\al})$. Besides,
by our choice of $\al$, we have $\tal(\tc_+(\R)) = \tc_+(\R)$.
Since $\tal$ preserves the orientation, the claim follows.

\vspace{2mm}

\noindent c) Let $p\in M\setminus Z_{\al}$ be an arbitrary point.
The homology class of $\alpha$ is nontrivial if and only if the
set $M\setminus Z_{\alpha}$ is connected. Therefore, there is a
regular $C^1$-path $\ga:[0,1]\to M$ such that $\ga(0)=c_+(0)$,
$\ga((0,1))\subset M\setminus Z_{\al}$, $\ga(1)=p$ and the pair
$(\dot{c}_+(0), \dot{\ga}(0))$ is positively oriented.

Let $\tga:[0,1]\to D$ be the lift of $\gamma$ to $D$ starting at
$\tc_+(0)$. Then $\tga((0,1])\subset E$, and hence $\tga(1)\in E$.
Since $\pi(\tga(1))=p$, the claim follows.
\end{proof}
\end{lem}

\medskip

In the following lemma we crucially use that the genus of $M$ is
greater than one.

\begin{lem}              \label{simple_curve_lem}
For every pair $(p,q)\in M\times M$ there exists a free homotopy class
$\alpha$ of closed curves in $M$ with the following properties:
\begin{itemize}
\item[a)]
The class $\al$ is homologically nontrivial and contains a simple
closed curve;
\item[b)]
We have $\{p,q\}\subset M\setminus Z_{\al}$.
\end{itemize}
\begin{proof}
Since the genus of $M$ is at least two, there exist free homotopy
classes  of closed curves in $M$, say $\al_1, \al_2, \al_3$, that
satisfy condition a) and have pairwise disjoint representatives.
Moreover, we can assume that the homology classes of $\al_1,
\al_2, \al_3$ are pairwise linearly independent. See
figure~~\ref{genus_two}.

%
%\begin{figure}[htbp]
%\begin{center}
%\input{gen_mor_two.pstex_t}
%\caption{Disjoint closed curves in a surface of genus at least two
%whose homology classes satisfy our conditions.} \label{genus_two}
%\end{center}
%\end{figure}
%

Let $Z_{\al_1},\, Z_{\al_2}$ and $Z_{\al_3}$ be the corresponding
subsets of $M$. See the discussion following
Theorem~~\ref{two_minim_thm}. It follows from
Theorem~~\ref{no_intersect_thm}  that the sets $Z_{\al_1},\,
Z_{\al_2}$ and $Z_{\al_3}$ are pairwise disjoint. Therefore at
least one of the classes $\al_1, \al_2, \al_3$ satisfies condition
b).
\end{proof}
\end{lem}

\medskip

We will now prove the main result of this section.

\begin{thm}              \label{high_genus_thm}
Every closed riemannian surface $M$ of genus greater than one is
totally insecure.
\begin{proof}
Let $(p,q)\in M\times M$ be an arbitrary
pair of points. We choose a free homotopy class $\alpha$ according
to Lemma~~\ref{simple_curve_lem};  let $E\subset D$ be as in
Lemma~~\ref{convex_lem}. By claim c) of Lemma~~\ref{convex_lem},
there exist points $x,y\in E$ such that $\pi(x)=p$, $\pi(y)=q$.
Let $\tal:D\to D$ be the covering transformation from
Lemma~~\ref{convex_lem} b). For $n\in\Z$ we set $y_n=\tal^n(y)$.

\medskip

By claim b) of Lemma~~\ref{convex_lem},  $y_n\in E$ for all
$n\in\Z$. For every positive integer $n$ let $\tc_n:[0,L_n]\to D$
be a $\tg$-minimal geodesic from $x$ to $y_n$. By claim a) of
Lemma~~\ref{convex_lem}, all curves $\tc_n([0,L_n])$ are contained
in $E$. We will show that the projections $c_n=\pi\circ\tc_n:[0,
L_n]\to M$ satisfy conditions a)-d) of
Proposition~~\ref{suffice_prop}, if we set $A=Z_{\alpha}$.

\medskip

Let $\td(\cdot,\cdot)$) be the distance in $D$ corresponding to
the metric $\tg$. Since
$\lim\limits_{n\to\infty}\td(x,y_n)=\infty$, we have $\lim
L_n=\infty$, i. e.~~condition a) of
Proposition~~\ref{suffice_prop} is satisfied. The curves
$\tc_n([0,L_n])$ are contained in $E$ and, by claim c) of
Lemma~~\ref{convex_lem}, $\pi(E)\subset M\setminus A$; hence
condition b) holds. Since the geodesics $\tc_n$ are $\tg$-minimal,
the geodesic segments $c_n([0, L_n))$ have no conjugate points,
which verifies condition c).

\medskip

Suppose now that condition d) of Proposition~~\ref{suffice_prop}
is not fulfilled. Then there is $\ve>0$, such that for all
$k\in{\mathbb N}$ there is an integer $n_k>0$ and a number
$t_k\in[k, L_{n_k}-k]$ such that
\begin{equation}       \label{low_bound_eq}
\td(\tc_{n_k}(t_k),\tc_+(\R)) \ge \ve.
\end{equation}
Observe that $\td(\tc_n(L_n),\tc_+(\R))=\td(y,\tc_+(\R))<\infty$.
Theorem~~\ref{bound_dist_thm} and standard facts from hyperbolic
geometry imply that there exists a constant $R$ such that
for all $n\in\N$ and all $t\in[0, L_n]$ we have
\begin{equation}       \label{upper_bound_eq}
\td(\tc_n(t),\tc_+(\R)) \le R.
\end{equation}

\medskip

\begin{figure}[htbp]
\begin{center}
\input{gen_mor_two.pstex_t}
\caption{Illustration to the proof of
Lemma~~\ref{simple_curve_lem}.} \label{genus_two}
\end{center}
\end{figure}

The group $\{\tal^n:n\in\Z\}$ acts cocompactly on the set $\{z\in
D:\,\td(z,\tc_+(\R))\le R\}$. Therefore, there is a sequence
$m_k$ such that the sequence of geodesics
$$
t\to\tal^{m_k}\circ\tc_{n_k}(t_k+t)
$$
converges to a $\tg$-minimal geodesic $c:\R\to\bar{E}$.

\medskip

By equation~~\eqref{upper_bound_eq}, the geodesic $c$ satisfies
$\td(c(t),\tc_+(\R))\le R$ for all $t\in\R$; hence
\begin{equation}      \label{endpts_eq}
\{c(\infty),c(-\infty)\} = \{\tc_+(\infty),\tc_+(-\infty)\}.
\end{equation}

\medskip

Equation~~\eqref{endpts_eq} and Theorem~~\ref{two_minim_thm} imply
that the trace of $c$ is either equal to the trace of $\tc_+$ (in
Case (i)) or is contained in the closed strip $S$ bounded by the
traces of $\tc_+$ and $\tc_-$ (in Case (ii)). Since
$c(\R)\subset\bar{E}$, we conclude that $c=\tc_+$, up to
parameterization. On the other hand, since $c(0)$ is the limit point of the
sequence $\tal^{m_k}\circ\tc_{n_k}(t_k)$, by
equation~~\eqref{low_bound_eq}, we have
$\td(c(0)),\tc_+(\R))\ge\ve>0$, a contradiction.

\medskip

We have shown that condition d) of Proposition~~\ref{suffice_prop}
also holds. Thus, Proposition~~\ref{suffice_prop} applies and
establishes the claim.
\end{proof}
\end{thm}


\begin{thebibliography}{99}
%
%\bibitem[Ar]{Ar} V.I. Arnold, {\it Mathematical methods of
%classical mechanics}, Nauka, Moscow,  1979 (the Russian original).
%\medskip
\bibitem[Ba1]{Ba1}
V.~Bangert, {\it Mather sets for twist maps and geodesics on
tori},  Dynamics Reported {\bf 1} (1988), 1 -- 56.
%Chichester-Stutt\-gart: John Wiley und
%B.G. Teubner 1988.
\medskip
\bibitem[Ba2]{Ba2}
V.~Bangert, {\it Geodesic rays, Busemann functions and monotone
twist maps}, Calc. Var.~{\bf 2} (1994), 49 -- 63.
\medskip
\bibitem[Ba3]{Ba3} V.~Bangert, {\em Hypersufaces without
self-intersections in the torus}, pp. 55 -- 71 in {\em Twist
Mappings and their Applications}, IMA Proceedings, vol. 44,
Springer-Verlag, Berlin, 1992.
\medskip
\bibitem[Ba4]{Ba4}
V.~Bangert, {\em On riemannian two-tori with a pole}, manuscript in preparation
(2009).
\medskip
\bibitem[BaG08]{BaG08} V. Bangert and E. Gutkin, {\em Secure two-dimensional tori are flat},
preprint arXiv:0806.3572 (2008).
\medskip
\bibitem[BG08]{BG08} K. Burns and E. Gutkin, {\em Growth of the number of geodesics between
points and insecurity for riemannian manifolds}, Discr. \& Cont. Dyn. Sys. {\bf A 21} (2008), 403 -- 413.
\medskip
\bibitem[Bu]{Bu}
H.~Busemann, {\it The Geometry of Geodesics}, Academic
Press, New York, 1955.
\medskip
\bibitem[GKM]{GKM} D. Gromoll, W. Klingenberg, W. Meyer, {\em
Riemannsche Geometrie im Grossen}, Lecture Notes in Mathematics
{\bf 55}, Springer-Verlag, Berlin, 1968.
\medskip
\bibitem[Gut05]{Gut05} E. Gutkin, {\em Blocking of billiard orbits and security for polygons and flat surfaces},
GAFA: Geom. \& Funct. Anal. {\bf 15} (2005), 83 -- 105.
\medskip
%\bibitem[Gut06]{Gut06} E. Gutkin, {\em Insecure configurations in lattice translation surfaces,
%with applications to polygonal billiards}, Discr. \& Cont. Dyn. Sys.
%{\bf A 16} (2006), 367 -- 382.
\bibitem[Gut09]{Gut09} E. Gutkin, {\em Topological entropy and blocking cost for
geodesics in Riemannian manifolds}, Geom. Dedicata {\bf 138}
(2009), 13 -- 23.
\medskip


%\bibitem{GJ} E. Gutkin and C. Judge, {\em Affine mappings of translation surfaces:
%geometry and  arithmetic}, Duke Math. J. {\bf 103} (2000), 191 --
%213.
%\bibitem {GK} E. Gutkin and A. Katok, {\em Caustics for inner and outer billiards},
%Comm. Math. Phys. {\bf 173} (1995), 101 -- 133.

\medskip
\bibitem[GS06]{GS06} E. Gutkin and V. Schroeder, {\em Connecting geodesics and
security of configurations in compact locally symmetric spaces},
Geom. Dedicata {\bf 118} (2006), 185 -- 208.
\medskip
\bibitem[He]{He}
G.A.~Hedlund, {\it Geodesics on a two-dimensional Riemannian
manifold with periodic coefficients}, Ann. Math.~{\bf 33} (1932),
719 -- 739.
\medskip
\bibitem[Her]{Her}
P.~Herreros, {\it Blocking: new examples and properties of products}, Erg. Theory \& Dyn. Sys.~{\bf 29} 
(2009), 569 -- 578.
\medskip
\bibitem[Hi]{Hi} M.W. Hirsch, {\em Differential Topology},
Springer-Verlag, New York, 1976.
\medskip
\bibitem[Ho]{Ho}
E.~Hopf, {\it Closed surfaces without conjugate points}, Proc.
Nat.~Acad.~Sci.~USA {\bf 34} (1948), 47 -- 51.
\medskip
\bibitem[In]{In}
N.~Innami, {\it Families of geodesics which distinguish flat
tori}, Math. J. Okayama Univ.~{\bf 28} (1986), 207 -- 217.
\medskip
\bibitem[Kl]{Kl} W. Klingenberg, {\it A course in differential
geometry}, Springer-Verlag, New York, 1978.
\medskip
\bibitem[LS07]{LS07} J.-F. Lafont and B. Schmidt, {\em Blocking light in compact Riemannian manifolds},
Geometry \& Topology {\bf 11} (2007), 867 -- 887.
\medskip
%\bibitem[Me97]{Mane} R. Ma\~n\'e, {\em On the topological entropy of geodesic flows},
%J. Diff. Geom.  {\bf 45} (1997), 74 -- 93.
%\medskip
%\bibitem[Ma79]{Manning} A. Manning, {\em Topological entropy for geodesic flows}, Ann. Math.
%{\bf 110}  (1979), 567 -- 573.
%\medskip
\bibitem[Mat]{Mat} J.N. Mather, {\em Action minimizing invariant
measures for positive definite lagrangian systems}, Math. Z. {\bf
207} (1991), 169 -- 207.
\medskip
\bibitem[Mo]{Mo}
M.~Morse, {\it A fundamental class of geodesics on any closed
surface of genus greater than one}, Trans. Amer.~Math.~Soc.~{\bf
26} (1924), 25 -- 60.
\medskip
\bibitem[Sch]{Sch} S. Schwartzman, {\em Asymptotic cycles}, Ann.
Math. {\bf 66} (1957), 270 -- 284.
\medskip
\bibitem[So]{So} B. Solomon, {\em On foliations of $\R^{n+1}$ by
minimal hypersurfaces}, Comment. Math. Helv. {\bf 61} (1986), 67
-- 83.
\medskip
\bibitem[We]{We} A. Weil, {\em On systems of curves on a
ring-shaped surface}, J. Indian Math. Soc. {\bf 19} (1931), 109
-- 114.

\bibitem[WKH]{WKH} Wing Kai Ho, {\it On Blocking Numbers of
Surfaces}, preprint arXiv:0807.2934 (2008).
\medskip
\end{thebibliography}
\end {document}